\documentclass[11pt,twoside]{article}
\usepackage{latexsym}
\usepackage{amssymb,amsbsy,amsmath,amsfonts,amssymb,amscd,mathtools}
\usepackage{graphicx}
\usepackage{hyperref}
\usepackage{pstricks}
\newcommand{\commentout}[1]{}
\newcommand{\R}{\mathbb{R}}

\newcommand {\e}  {\varepsilon}

\newcommand {\da} {\delta}

\newcommand {\vp} {\varphi}

\newcommand {\Chi} {{\bf \raise 2pt \hbox{$\chi$}} }

\newcommand {\f}   {\frac}
\newcommand {\p}   {\partial}
\newcommand{\fer}{\eqref}

\newcommand {\proof} {\noindent {\bf Proof}. }
\newcommand{\beq}{\begin{equation}}
\newcommand{\eeq}{\end{equation}}
\newcommand{\bea} {\begin{array}{rl}}
\newcommand{\eea} {\end{array}}
\newcommand{\bepa}{\left\{ \begin{array}{l}}
\newcommand{\eepa} {\end{array}\right.}
\newtheorem{theorem}{Theorem}[section]
\newtheorem{lemma}[theorem]{Lemma}

\newtheorem{prop}[theorem]{Proposition}
\newtheorem{corollary}[theorem]{Corollary}
\newcommand{\qed}{{ \hfill
                       {\unskip\kern 6pt\penalty 500 \raise -2pt\hbox{\vrule\vbox to 6pt{\hrule width 6pt
                       \vfill\hrule}\vrule} \par}   }}

\title{A Hamilton-Jacobi approach to characterize the evolutionary equilibria in heterogeneous environments}

\author{
Sepideh Mirrahimi\thanks{Institut de Math\'ematiques de Toulouse; UMR 5219, Universit\'e de Toulouse; CNRS, UPS IMT, F-31062 Toulouse Cedex 9, France; E-mail: Sepideh.Mirrahimi@math.univ-toulouse.fr} }
 
\date{\today}

\begin{document}
\maketitle
\pagestyle{plain}
\pagenumbering{arabic}

\begin{abstract}
In this work, we characterize the solution of a system of   elliptic integro-differential equations describing a phenotypically structured population subject to mutation, selection and migration between two habitats. Assuming that the effects of the mutations are small but nonzero, we show that the population's distribution has at most two peaks and we give explicit conditions under which the population will be monomorphic (unimodal distribution) or dimorphic (bimodal distribution). More importantly, we provide a general method to determine the dominant terms of the population's distribution in each case. Our work, which is based on Hamilton-Jacobi equations with constraint, goes further than previous works where such tools 
were used, for different problems from evolutionary biology, to identify the asymptotic solutions, while the mutations vanish, as a sum of Dirac masses. In order to extend such results to the case with non-vanishing effects of mutations, the main elements are a uniqueness property and the computation of the correctors.\\
This method allows indeed to go further than the Gaussian approximation commonly used by biologists and makes a connection between the theories of adaptive dynamics and
quantitative genetics. Our work being motivated by biological questions, the objective of this article is to provide the mathematical details which are necessary for our biological results \cite{SG.SM:17}.
\end{abstract}

\bigskip

\section{Introduction}

\noindent 
Can we characterize the phenotypical distribution of a population which is subject to the Darwinian evolution? The mathematical modeling of the phenotypically structured populations, under the effects of mutations and selection leads to parabolic and elliptic integro-differential equations. The solutions of such equations, as the mutation term vanishes, converge to a sum of Dirac masses, corresponding to the dominant traits. During the last decade, an approach based on Hamilton-Jacobi equations with constraint has been developed which allows to describe such asymptotic solutions.  There is a large literature on this method. We refer to \cite{OD.PJ.SM.BP:05,GB.BP:08,SM:11} for the establishment of the basis of this approach for problems from evolutionary biology. Note that related tools were already used in the case of local equations (for instance KPP type equations) to describe the propagation phenomena (see for instance \cite{MF:86,LE.PS:89}). 

\noindent
Such results, which are based on a logarithmic transformation (the so-called Hopf-Cole transformation) of the population's density, provide mainly the convergence along subsequences of the logarithmic transform to a viscosity solution of a Hamilton-Jacobi equation with constraint, as the effects of the mutations vanish. This allows to obtain a qualitative description of the population's phenotypical distribution for vanishing mutations' steps. To be able to characterize the population's distribution for non-vanishing effects of mutations, one should prove a uniqueness property for the viscosity solution of the Hamilton-Jacobi equation with constraint and compute the next order terms. Such properties are usually not studied due to technical difficulties. However, from the biological point of view it is usually more relevant to consider non-vanishing mutations' steps. 

\noindent
In this work, as announced in \cite{SG.SM:16}, we provide such analysis, including a uniqueness result and the computation of the correctors, in the case of a selection, mutation and migration model.  Note that a recent work \cite{SM.JR:15-1,SM.JR:16} has also provided similar results in the case of homogeneous environments. We believe indeed that going further in the Hamilton-Jacobi approach for different problems from evolutionary biology, by providing higher order approximations, can make this approach more useful for the evolutionary biologists.

\noindent
 The purpose of this article is to provide the mathematical details and proofs which are necessary for our biological results \cite{SG.SM:17}.  
As explained in \cite{SG.SM:17}, our method allows to  provide more quantitative results and correct the previous approximations obtained by biologists. 
\\

\noindent
Our objective is to characterize the  solutions to the following system, for $z\in  \R$,
\begin{equation}
\label{main}
\begin{cases}
-\e^2  n_{\e,1}''(z)=n_{\e,1}(z) R_1(z,N_{\e,1})+m_2n_{\e,2}(z)-m_1n_{\e,1}(z),\\
-\e^2 n_{\e,2}''(z)=n_{\e,2}(z) R_2(z,N_{\e,2})+m_1n_{\e,1}(z)-m_2n_{\e,2}(z),\\
N_{\e,i}=\int_\R n_{\e,i}(z)dz, \qquad \text{for $i=1,2$},
\end{cases}
\end{equation}
with 
\beq
\label{Ri}
R_i(z,N_i)= r_i- g_i(z-\theta_i)^2-\kappa_i N_i,\qquad \text{with $\theta_1=-\theta$ and $\theta_2=\theta$}.
\eeq
This system represents the equilibrium of a population that is structured by a phenotypical trait $z$, and which is subject to selection, mutation and migration between two habitats. We denote by $n_i(z)$ the density of the phenotypical distribution in habitat $i$, and by $N_i$ the total population size in habitat $i$. The growth rate $R_i(z,N_i)$ is given by \fer{Ri}, where $r_i$ represents the maximum intrinsic growth rate, the positive constant $g_i$ is the strength of the selection, $\theta_i$ is the optimal trait in habitat $i$ and the positive constant $\kappa_i$ represents the intensity of the competition. The nonnegative constants $m_i$ are the migration rates between the habitats. 

\noindent
Such phenomena have already  been studied using several approaches by the theoretical evolutionary biologists. A first class of results  are based on the adaptive dynamics approach, where one considers that the mutations are very rare such that the population has time to attain its equilibrium between two mutations and hence the population's distribution has discrete support (one or two points in a two habitats model) \cite{GM.IC.SG:97,TD:00,CF.SM.EM:12}.  A second class of results are based on an approach known as  'quantitative genetics', which allows more frequent mutations and does not  separate the evolutionary and the ecological time scales so that the population's distribution is continuous (see \cite{Rice-book}--chapter $7$). A main assumption in this class of works is that one considers that the population's distribution is a gaussian \cite{AH.TD.ET:01,OR.MK:01} or, to take into account the possibility of dimorphic populations, a sum of one or two gaussian distributions \cite{SY.FG:09,FD.OR.SG:13}. 

\noindent
In our work, as in the  quantitative genetics framework, we also consider continuous phenotypical distributions. However,  we don't assume any a priori gaussian assumption. We compute directly the population's distribution and in this way we correct the previous approximations.  To this end, we also provide some results in the framework of adaptive dynamics and in particular, we generalize previous results on the identification of the evolutionary stable strategy (ESS) (see Section \ref{sec:ad} for the definition) to the case of nonsymetric habitats. Furthermore, our work makes a connection between the two approaches of adaptive dynamics and quantitative genetics.
\\

\noindent
{\bf Assumptions:}\\
To guarantee that the population does not get extinct, we assume that
\begin{equation}
\label{as:r-m}
\max( r_1-m_1,r_2-m_2) > 0.
\end{equation}
 Moreover, in the first part of this article, we  assume that there is positive migration rate in both directions, i.e.
\begin{equation}
\label{as:m}
m_i>0, \qquad \text{i=1,2}.
\end{equation}
The source and sink case, where for instance $m_2=0$, will be analyzed in the last section. \\

\noindent 
Note that  in \cite{SM:12} the limit, as $\e\to 0$ and along subsequences, of the solutions to such system, under assumption \fer{as:m}, and in a bounded domain, was studied. In the present work, we go further than the asymptotic limit along subsequences and we obtain uniqueness of the limit and identify the dominant terms of the solution when $\e$ is small but nonzero.  In this way, we are able to characterize the solution when the mutation's steps are not negligible. 
\\

\noindent
 {\bf The main elements of the method:}\\
To describe the solutions $n_{\e,i}(z)$ we use a WKB ansatz
\beq
\label{WKB}
n_{\e,i}(z)=\frac{1}{\sqrt{2\pi \e}} \exp \left(\frac{u_{\e,i}(z)}{\e} \right).
\eeq
Note that a first approximation that is commonly used in the theory of 'quantitative genetics', is a gaussian distribution of the following form
$$
n_{\e,i}(z)=\frac{N_i}{\sqrt{2\pi\e}\sigma}\exp  \left(\frac{-(z-z^*)^2}{\e\sigma^2} \right)
=\frac{1}{\sqrt{2\pi\e}} \exp \left(\frac{-\frac{1}{2\sigma^2}(z-z^*)^2+\e\log \frac{N_i}{\sigma}}{\e} \right).
$$
Here, we try to go further than this a priori gaussian assumption and to approximate directly $u_{\e,i}$. To this end, we write an expansion for 
$u_{\e,i}$ in terms of $\e$:
\beq
\label{ap-ue}
u_{\e,i}= u_i+\e v_i+\e^2 w_i+O(\e^3).
\eeq
We first prove that $u_1=u_2=u$ is the unique viscosity solution to a Hamilton-Jacobi equation with constraint which can be computed explicitly. The uniqueness of solution of such Hamilton-Jacobi equation with constraint is related to the uniqueness of the ESS  and  to the weak KAM theory \cite{AF:16}.
Such function  $u$ indeed satisfies 
$$
\max_\R \; u(z)= 0,
$$
with the maximum points attained at one or two points corresponding to the ESS points of the problem.
We then notice that, while $u(z)<0$, $n_{\e,i}(z)$ is exponentially small. Therefore, only the values of $v_i$ and $w_i$ at the points which are close to the zero level set of $u$ matter, i.e. the ESS points. We next show how to compute formally $v_i$ and hence its second order Taylor expansion around the ESS points,  and the value of $w_i$ at those points. These approximations together with  a fourth order Taylor expansion of $u_i$ around the ESS points are indeed enough to approximate the moments of the population's distribution with an error of order $\e^2$.
\\

\noindent
The paper is organized as follows. In Section \ref{sec:ad} we introduce some notions from the theory of adaptive dynamics that will be used in the following sections. In Section \ref{sec:results} we state our main results (theorems \ref{th:ESS} and \ref{thm:main}) and discuss their consequences.  In this section, we also provide the method to compute the correctors and approximate the moments of the population's distribution. In Section \ref{sec:pr-ESS} we provide the proofs of the results in the adaptive dynamics framework and in particular we prove Theorem \ref{th:ESS}. In Section \ref{sec:pr-HJ} we prove Theorem \ref{thm:main}. Finally, in Section \ref{sec:sink} we generalize our results to the sink and source case where the migration is only in one direction ($m_2=0$).

\section{Some notions from the theory of adaptive dynamics }
\label{sec:ad}

In this section, we introduce some notions  from the theory of adaptive dynamics  that we will be using in the next sections  \cite{GM.IC.SG:97}.   Note that our objective is not to study the framework of adaptive dynamics where the mutations are assumed to be very rare. However, these notions appear naturally from our asymptotic computations.\\

\noindent
{\bf Effective fitness:} 
 The effective fitness $W(z;N_1,N_2)$ is the largest eigenvalue of the following matrix:
\beq
\label{efitness}
\mathcal A (z;N_1,N_2)= \left( 
\begin{array}{cc}
R_1(z ; N_1) -m_1 & m_2\\
m_1 & R_2(z ; N_2) -m_2
\end{array}
\right),
\eeq
that is
\beq
\label{W}
\begin{array}{rl}
W(z;N_1,N_2) &= \f 1 2 \Big[ \left( R_1(z ; N_1)+R_2(z ; N_2)-m_1 - m_2 \right) \\
 &+\sqrt{ \left( R_1(z ; N_1)-R_2(z ; N_2)-m_1 + m_2 \right)^2
+4m_1m_2 ) } \, \Big].
\end{array}
\eeq
This indeed corresponds to the \emph{effective} growth rate associated with trait $z$ in the whole metapopulation when the total population sizes are given by $(N_1,N_2)$.
\\

 \noindent
{\bf Demographic equilibrium:} 
Consider a set of points $\Omega=\{z_1,\cdots z_m\}$. The demographic equilibrium corresponding to this set is given by $(n_1(z),n_2(z))$, with the total population  sizes $(N_1,N_2)$, such that
$$
n_i(z)=\sum_{j=1}^m \alpha_{i,j}\delta(z-z_j),\quad N_i = \sum_{j=1}^m \alpha_{i,j},\quad W(z_j,N_1,N_2)=0,\quad i=1,2,\; j=1,\cdots, m,
$$
and such that $(\alpha_{1,j},\alpha_{2,j})^{T}$  is the right eigenvector associated with the largest eigenvalue $W(z_j,N_1,N_2)=0$ of $\mathcal A(z_j;N_1,N_2)$.
\\

\noindent 
{\bf Invasibility:} We say that a mutant trait $z_m$ can invade a resident strategy $\{ z^M \}$ at its demographic equilibrium $( N_1^M, N_2^M)$ if $W(z_m, N_1^M,N_2^M)>0$.
 \\
 
\noindent
{\bf Evolutionary stable strategy:} 
A set of points $\Omega^*=\{z_1^*,\cdots, z_m^*\}$ is called an evolutionary stable strategy (ESS) if
$$
W(z,N^*_1,N_2^*)=0,\quad \text{for $z\in \mathcal A$ and} ,\quad W(z,N_1^*,N_2^*)\leq 0, \quad \text{for  $z\not\in \mathcal A$,}
$$
where $N_1^*$ and $N_2^*$ are the total population sizes corresponding to the demographic equilibrium associated with the set $\Omega^*$.
\\

\noindent 
{\bf Notation}: We will use the star sign $^*$ whenever we talk about an evolutionary stable strategy $\Omega^*$ (and similarly for the corresponding demographic equilibrium $(n_1^*,n_2^*)$ and the total population sizes $(N_1^*,N_2^*)$). We add an index $M$  when the strategy is monomorphic (a set of a single trait $\{z^{M*}\}$ with the corresponding demographic equilibrium $(n_1^{M*},n_2^{M*})$, and the total population sizes $(N_1^{M*}, N_2^{M*})$) and an index $D$ when the strategy is dimorphic (a set of two traits $\{z_{\mathrm{I}}^{D*},z_{\mathrm{II}}^{D*}\}$ with the corresponding demographic equilibrium $(n_1^{D*}, n_2^{D*})$, and the total population sizes $(N_1^{D*},N_2^{D*})$).

\section{The main results and the details of the method}
\label{sec:results}

In this section, we state our main results and provide the details of our method for the approximation of the equilibrium distribution $n_{\e,i}(z)$. In Subsection \ref{sec:ad-results} we provide the results in the framework of adaptive dynamics. In Subsection \ref{sec:u} we state our main result on the convergence to the zero order term $u_i$ and its explicit computation. In Subsection \ref{sec:vw} we show how to compute the next order terms. Finally, in Subsection \ref{sec:mom} we provide the approximation of the moments of the population's distribution.

\subsection{The adaptive dynamics framework}
\label{sec:ad-results}

Our main result in the adaptive dynamics framework is that there exists a unique ESS which is whether monomorphic (a single Dirac mass) or dimorphic (a sum of two Dirac masses).  
  We determine indeed under which conditions the ESS is monomorphic or dimorphic. To state our result, we first define
\beq
\label{dim-eq}
 z^{D*}=\sqrt{\theta^2-\f{m_1m_2}{4\theta^2g_1g_2}},\qquad 
 N_1^{D*}=\f{ \f{m_1m_2}{4\theta^2g_2}+r_{2}-m_1}{\kappa_1},\qquad   N_2^{D*}=\f{ \f{m_1m_2}{4\theta^2g_1}+r_{2}-m_2}{\kappa_2}.
\eeq

\begin{theorem}
\label{th:ESS}
Assume  \fer{as:r-m}--\fer{as:m}. Then, there exists a unique set of points $\Omega^*$ which is an ESS.  \\
(i) The ESS is dimorphic if and only if
\beq
\label{as1:dim}
\f{m_1m_2}{4g_1g_2 \theta^4}< 1,
\eeq
\beq 
\label{as2:dim}
0<m_2  N_2^{D*} +(R_1(- z^{D*}; N_1^{D*})- m_1)   N_1^{D*},
\eeq
and
\beq 
\label{as3:dim}
0<m_1   N_1^{D*} +(R_2( z^{D*};N_2^{D*})- m_2)  N_2^{D*}.
\eeq
Then the dimorphic  equilibrium is given by
\beq
\label{dim-ESS}
n_i^{D*}=\nu_{\mathrm{I},i} \da(z+ z^{D*})+ \nu_{\mathrm{II},i} \da(z- z^{D*}),\quad  \nu_{\mathrm{I},i}+\nu_{\mathrm{II},i}= N_i^{D*},\quad i=1,2.
\eeq
(ii) If the above conditions are not satisfied then the ESS is monomorphic. In the case where  condition \fer{as1:dim}  is verified but the r.h.s. of  \fer{as2:dim} (respectively  \fer{as3:dim}) is negative,  the fittest trait belongs to the interval $(- \theta, - z^{D*})$ (respectively  $( z^{D*}, \theta)$). If  \fer{as1:dim} is satisfied but  \fer{as2:dim} (respectively  \fer{as3:dim}) is an equality then the monomorphic ESS is given by $\{-z^{D*}\}$ (respectively $\{z^{D*}\}$).

\end{theorem}

\medskip
 
 \noindent
Note that one can compute the weights $\nu_{k,i}$, for $k=\mathrm{I},\mathrm{II}$ and $i=1,2$:
\beq
\label{nuij}
\begin{array}{c}
\left(
\begin{array}{c}
\nu_{\mathrm{I},1}\\
\nu_{\mathrm{I},2} 
\end{array}
\right)
= \f{ m_1   N_1^{D*} +(R_2( z^{D*};N_2^{D*})- m_2)  N_2^{D*}}{ m_1m_2 - \big(R_1(- z^{D*};  N_1^{D*})-m_1 \big) \big(R_2( z^{D*}; N_2^{D*})-m_2 \big)  }\left(
\begin{array}{c}
m_2\\
-R_1(- z^{D*} ;  N_1^{D*})+m_1
\end{array}
\right),
\\
\left(
\begin{array}{c}
\nu_{\mathrm{II},1} \\
\nu_{\mathrm{II},2} 
\end{array}
\right)
=\f{m_2  N_2^{D*} +(R_1(- z^{D*}; N_1^{D*})- m_1)   N_1^{D*}}{  m_1m_2 - \big(R_1(- z^{D*};  N_1^{D*})-m_1 \big) \big(R_2( z^{D*}; N_2^{D*})-m_2 \big) }   
\left(
\begin{array}{c}
-R_2( z^{D*} ; N_2^{D*})+m_2\\
m_1
\end{array}
\right).
\end{array}
\eeq

\noindent
Moreover,  since $W(- z^{D*};  N_1^{D*},N_2^{D*} )=0$, one can easily verify that 
 condition \fer{as2:dim} is equivalent with
\beq
\label{as:eq7}
m_1  N_1^{D*} +  (R_2(- z^{D*};N_2^{D*})- m_2)  N_2^{D*}<0.
\eeq
Similarly, since $W( z^{D*};  N_1^{D*},N_2^{D*} )=0$, one can easily verify that 
 condition \fer{as3:dim} is equivalent with
\beq
\label{as:eq8}
 m_2 N_2^{D*} +  (R_1( z^{D*}; N_1^{D*})- m_1)   N_1^{D*}<0.
\eeq
 
 \noindent
To prove Theorem  \ref{thm:main}--(iii) we will use the following result which is a corollary of Theorem \ref{th:ESS}.

\begin{corollary}
\label{cor:deg} 
Assume that
\beq
\label{non-deg}
m_2  N_2^{D*} +(R_1(- z^{D*}; N_1^{D*})- m_1)   N_1^{D*}\neq 0,
\qquad
m_1   N_1^{D*} +(R_2( z^{D*};N_2^{D*})- m_2)  N_2^{D*}\neq 0,
\eeq
and let the set $\Omega^*$ be the unique ESS of the model and $(N_1^*,N_2^*)$ be the total population sizes at the demographic equilibrium of this ESS. Then, 
\beq
\label{Wneg}
W(z,N_1^*,N_2^*)<0,\qquad \text{for all z $\in \R \setminus \Omega^*$}.
\eeq

\end{corollary}

\noindent
Note also  that when the habitats are symmetric, then conditions \fer{as2:dim} and \fer{as3:dim} always hold under condition \fer{as1:dim}, and hence
\begin{corollary} \label{cor-sym}
Assume that the habitats are symmetric:
\beq
\label{par-sym}
r=r_{1}=r_{2},\quad g=g_1=g_2, \quad \kappa=\kappa_1=\kappa_2,\quad m=m_1=m_2.
\eeq
(i) Then the unique ESS is dimorphic if and only if
\beq
\label{dim-sym}
\f{m }{2g}< \theta^2.
\eeq
 The dimorphic ESS is determined by \fer{dim-ESS}.\\
(ii) When condition \fer{dim-sym} is not satisfied, then the ESS is monomorphic and the corresponding monomorphic equilibrium is given by
\beq
\label{ESS-sym-mono}
 n_1^{M*}(z)= n_2^{M*}(z)=  N^{M*} \, \da(z), \quad \text{with } N^{M*}=\f{1}{\kappa} \left( r -g  \theta^2\right).
\eeq

\end{corollary}

 \noindent
The next proposition gives an interpretation of conditions \fer{as2:dim} and \fer{as3:dim}.

\begin{prop}
\label{prop:invade}
Assume that condition \fer{as1:dim} is satisfied and that $r_i-m_i>0$, for $i=1,2$. Then,
\\
(i) condition \fer{as2:dim} holds if and only if a mutant trait of type $ z^{D*}$ can invade a monomorphic resident population of  type $- z^{D*}$ which is at it's demographic equilibrium.\\
(ii) condition \fer{as3:dim} holds if and only if a mutant trait of type $- z^{D*}$ can invade a  monomorphic resident population of  type $ z^{D*}$ which is at it's demographic equilibrium.
\end{prop}

\noindent
One can indeed rewrite conditions \fer{as2:dim} and \fer{as3:dim} respectively as below
$$
{C_1\, < \, \alpha_2 r_2 -  \alpha_1 r_1}
, \qquad 
{C_2\, < \, \beta_1 r_1 -\beta_2 r_2},
$$
with $C_i$, $\alpha_i$ and $\beta_i$ constants depending on $m_1$, $m_2$, $g_1$, $g_2$, $\kappa_1$, $\kappa_2$ and $\theta$. These conditions are indeed a measure of asymmetry between the habitats. They appear from the fact that even if condition \fer{as1:dim}, which is the only condition for dimorphism in symmetric habitats, is satisfied, while the quality of the habitats are very different,   the ESS cannot be dimorphic. In this case, the population will  be able to adapt only  to one of the habitats and it will be maladapted to the other one.

\subsection{The computation of the zero order terms $u_i$}
\label{sec:u}
The identification of the zero order terms $u_i$ is based on the following result.  

\begin{theorem}\label{thm:main}
Assume \fer{as:r-m}--\fer{as:m}. \\
(i) As $\e \to 0$, $(n_{\e,1}, n_{\e,2})$ converges to $(n_1^*,n_2^*)$, the demographic equilibrium of the unique  ESS of the model. Moreover, as $\e \to 0$, $N_{\e,i}$  converges to $N_i^*$, the total population size in patch $i$ corresponding to this demographic equilibrium.
\\
(ii) As $\e\to 0$, both sequences $(u_{\e,i})_\e$, for $i=1,2$, converge along subsequences and locally uniformly in $\R$ to a continuous function $u\in \mathrm{C}(\R)$,  such that $u$ is a viscosity solution to the following equation
\begin{equation}
\label{HJ}
\left\{ \begin{array}{ll}-|u'(z)|^2= W(z,N_1^*,N_2^*),&\quad \text{in $\R$},\\
\max_{z\in \R}u(z)=0.\end{array}\right.
\end{equation}
Moreover, we have the following condition on  the zero level set of $u$:
$$
{\rm supp}\, n_1^*= {\rm supp}\, n_2^*\subset \{z\, |\, u(z)=0 \} \subset \{z\, |\, W(z,N_1^*,N_2^*)=0\}.
$$
(iii) Under  condition \fer{non-deg} we have ${\rm supp}\, n_1^*= {\rm supp}\, n_2^*=  \{z\, |\, W(z,N_1^*,N_2^*)=0\}$ and hence
\beq
\label{aubry}
\{z\, |\, u(z)=0 \} = \{z\, |\, W(z,N_1^*,N_2^*)=0\}.
\eeq
The solution of \eqref{HJ}--\fer{aubry} is indeed unique and hence the whole sequence $(u_{\e,i})_\e$   converge  locally uniformly in $\R$ to   $u$.
\end{theorem}

\noindent
Note that a Hamilton-Jacobi equation of type \fer{HJ} in general might admit several viscosity solutions. Here, the uniqueness is obtained thanks to \fer{aubry} and a property from   the weak KAM theory, which is the fact that the viscosity solutions are completely determined by one value taken on each static class of the Aubry set (\cite{PL:82}, Chapter 5 and \cite{GC:01}).
 In what follows we assume that \fer{non-deg} and hence \fer{aubry} always hold. We then give an explicit formula for $u$ considering two cases (one can indeed verify easily that the functions below are viscosity solutions to  \eqref{HJ}--\fer{aubry}):\\
 
 \noindent
(i) {\bf Monomorphic ESS : } We consider the case where there exists a unique monomorphic ESS $z^{M*}$ and the corresponding demographic equilibrium is given by $(N_1^{M*}\da(z^{*}),N_2^{M*}\da(z^{M*}))$. Then  $u$ is given by
\beq
\label{u-exp}
u(z)= -\big| \int_{z^{M*}}^{ z} \sqrt{- W(x; N_1^{M*}, N_2^{M*})} dx \big|.
\eeq
(ii) {\bf Dimorphic ESS : } We next consider the case where there exists a unique dimorphic ESS $(z_{\mathrm{I}}^{D*},z_{\mathrm{II}}^{D*})$   with the demographic equilibrium:
$
n_i=\nu_{\mathrm{I},i} \da(z-z^{D*}_{\mathrm{I}})+\nu_{\mathrm{II},i} \da(z-z^{D*}_{\mathrm{II}}),
$ and
$ \nu_{\mathrm{I},i}+ \nu_{\mathrm{II},i} =N_i^{D*}$. Then  $u$ is given by
$$
 u(z)=\max \Big( - |\int_{z_{\mathrm{I}}^{*}}^{ z} \sqrt{- W(x; N_1^{D*}, N_2^{D*})} dx|
  , - |\int_{z_{\mathrm{II}}^{*}}^{z} \sqrt{- W(x; N_1^{D*}, N_2^{D*})} dx |\Big).
$$

\subsection{Next order terms}
\label{sec:vw}

%

\noindent
In this subsection we show how one can compute formally the first order term $v_i$, and in particular its second order Taylor expansion around  the zero level set of $u$, and determine the value of $w_i$ at those points. We only present the method in the case of monomorphic   population where the demographic equilibrium corresponding to this ESS  is given by $(N_1^{M*}\da(z-z^{M*}), N_2^{M*}\da(z-z^{M*}) )$. The dimorphic case can be treated following similar arguments.
\\
 
 \noindent
 We first note that, one can compute, using \fer{u-exp}, a Taylor expansion of order $4$ around the ESS point $z^{M*}$:
$$
u(z) = -\f{A}{2 }(z-z^{M*})^2+B(z-z^{M*})^3+C(z-z^{M*})^4+O(z-z^{M*})^5.
$$
 We then look for constants $D_i$, $E_i$, $F_i$ and $G_i$ such that
$$
v_i(z)=v_i(z^{M*})+D_i (z-z^{M*}) +E_i(z-z^{M*})^2+O(z-z^{M*})^3,\qquad w_i(z^{M*})=F_i+G_i(z-z^{M*})+O(z-z^{M*})^2.
$$
We will only compute $D_i$, $E_i$ and $F_i$. The constants $G_i$ are not necessary in the computation of the moments but they appear in our intermediate computations.
Replacing the functions $u$, $v_i$ and $w_i$ by the above approximations to compute $N_{\e,i}=\int_\R n_{\e,i} (z)dz$, we obtain
$$
v_i(z^{M*})= \log \big(N_i^{M*} \sqrt{A} \big),
$$
$$
N_{\e,i} =N_i^{M*} +\e K_i+O(\e^2),\quad \text{with }\quad K_i= N_i^{M*} \big(\f{7.5\, B^2}{A^3}+  \f{3(C+BD_i)}{A^2}+\f{E_i+0.5 \,D_i^2}{A}+F_i  \big).
$$

\noindent
Note also that writing \fer{main} in terms of $u_{\e,i}$ we obtain
\begin{equation}
\label{main2}
\begin{cases}
-\e  u_{\e,1}''(z) = |u_{\e,1}'|^2+ R_1(z,N_{\e,1})+m_2\exp \big( \f{u_{\e,2}-u_{\e,1}}{\e}\big)-m_1 ,\\
-\e  u_{\e,2}''(z) = |u_{\e,2}'|^2+ R_2(z,N_{\e,2})+m_1\exp \big( \f{u_{\e,1}-u_{\e,2}}{\e}\big)-  m_2.
\end{cases}
\end{equation}

\noindent
We  then let $\e\to 0$ in the first line of \fer{main2} and use \fer{HJ} to obtain 
\beq
\label{Q2/Q1}
  v_2(z)-v_1(z) = \log \Big( \f{1}{m_2}\big( W(z,N_1^{M*},N_2^{M*})-R_1(z,N_1^{M*})+m_1 \big) \Big).
\eeq
Keeping respectively, only the terms of order $(z-z^{M*})$ and $(z-z^{M*})^2$ we find
$$
\lambda_1=D_2-D_1=\f{2g_1N_1^{M*}(z^{M*}+\theta)} {m_2N_2^{M*}\, },
$$
$$
\lambda_2=E_2-E_1=\f{N_1^{M*}}{m_2N_2^{M*}}(-A^2+g_1) -\f{2g_1^2N_1^{M*\, 2}}{m_2^2N_2^{M* \, 2}} (z^{M*}+\theta)^2.
$$
Combining the above lines we obtain
\beq
\label{da1-da2}
\f{ K_2}{N_2^{M*} } -\f{ K_1}{N_1^{M*} } =\lambda_3+\f{0.5\,\lambda_1 (D_1+D_2)}{A}+F_2-F_1, \quad \text{with} \quad \lambda_3 =
\f{3B}{A^2}\lambda_1+
\f{1}{A}\lambda_2.
\eeq
 
\noindent
Next, keeping the terms of order $\e$ in \fer{main2} we obtain, for $\{i,j\} =\{1,2\}$,
\beq
\label{1-order}
-  u''=2  u' \cdot  v_i'-\kappa_i K_i +m_j \exp(v_j-v_i)(w_j-w_i).
\eeq
Evaluating the above equality at $z^{M*}$ we obtain
$$
A=-\kappa_iK_i+m_j \f{N_j^{M*}}{N_i^{M*}}(F_j-F_i).
$$
Replacing \fer{da1-da2} in the above system we obtain
$$
\begin{cases}
A=-\kappa_1 K_1+m_2 \f{N_2^{M*}}{N_1^{M*}}(\f{ K_2}{N_2^{M*} } -\f{ K_1}{N_1^{M*} } -\lambda_3-\f{0.5\,\lambda_1 (D_1+D_2)}{A} ),\\
A=-\kappa_2 K_2+m_1 \f{N_1^{M*}}{N_2^{M*}}(\f{ K_1}{N_1^{M*} } -\f{ K_2}{N_2^{M*} } +\lambda_3 +\f{0.5\,\lambda_1 (D_1+D_2)}{A}).
\end{cases}
$$
This system allows us to identify $(K_1,K_2)$ in a unique way, as an affine function of $(D_1+D_2)$.\\

\noindent
Next we substrate the two lines of the system \fer{1-order} to obtain 
\beq
\label{w2w1}
w_2-w_1= \f{2 u' \cdot  (v_2'-v_1')+\kappa_1 K_1-\kappa_2 K_2}{m_2\exp(v_2-v_1)+m_1\exp(v_1-v_2)}.
\eeq
Evaluating the above equation at $z^{M*}$ we find
$$
F_2-F_1= \f{\kappa_1 K_1-\kappa_2 K_2}{m_1N_1^{M*}/N_2^{M*} + m_2 N_2^{M*}/N_1^{M*}},
$$
and keeping the terms of order $(z-z^{M*})$ we obtain
$$
G_2-G_1=\f{-2A(D_2-D_1)}{m_1N_1^{M*}/N_2^{M*} + m_2 N_2^{M*}/N_1^{M*}}+\f{(m_2N_2^{M*}/N_1^{M*}-m_1N_1^{M*}/N_2^{M*})(D_2-D_1)}{(m_1N_1^{M*}/N_2^{M*} + m_2 N_2^{M*}/N_1^{M*})^2}\, (\kappa_1 K_1-\kappa_2 K_2).
$$
We then keep the terms of order $(z-z^{M*})$ in \fer{1-order} to find
$$
-6B=-2AD_1+m_2 \f{N_2}{N_1}\big((D_2-D_1)(F_2-F_1)+G_2-G_1\big).
$$
Combining the above lines, one can write $D_1$ as an affine function of $D_1+D_2$. Since $D_2-D_1$ is already known, this allows to identify, at least in a generic way, $D_i$ and consequently $K_i$ (see \cite{SG.SM:17} for examples of such computations).
Next, we replace \fer{w2w1} in \fer{1-order} to obtain
$$
-  u''=2  u' \cdot   v_i' - \kappa_i K_i+\f{m_j\exp(v_j-v_i)}{m_2\exp(v_2-v_1)+m_1\exp(v_1-v_2)} \big(2 u' \cdot  (v_j'-v_i')+\kappa_i K_i- \kappa_j K_j\big).
$$
All the terms in the above system, except  $v_i'$, are already known. Hence one can compute  $v_i$ from the above system. In particular, 
keeping the terms of order  $(z-z^{M*})^2$ in the above line, one can compute   $E_i=\f12v_i''(z^{M*})$ and consequently $F_i$.

\subsection{Approximation of the moments}
\label{sec:mom}

The above approximations of $u$, $v_i$ and $w_i$ around the ESS points allow us to estimate the moments of the population's distribution with an error of at most order $O(\e^2)$. We only provide such approximations in the monomorphic case. One can obtain such approximations in the case of dimorphic ESS following similar computations.  We first note that, replacing $u_{\e,i}$ by the approximation \fer{ap-ue} and using the Taylor expansions of $u$, $v_i$ and $w_i$ obtained above, we can compute
$$
\begin{array}{c}
\int (z-z^{M*})^k n_{\e,i}(z) dz = \f{\e^{\f k 2}\sqrt{A} N_i^{M*}}{\sqrt{2\pi}} \int_\R (y^k e^{-\f A 2 y^2} \big( 1+\sqrt{\e} (By^3+D_iy)+O(\e) \big) dy\\
=\e^{\f k 2} N_i^{M*}\Big( \mu_k(\f{1}{A})+\sqrt{\e}\big(B\mu_{k+3}(\f 1A)+D_i \mu_{k+1}(\f 1A) \big) +O(\e^{\f{k+2}{2}}),
\end{array}
$$
where $\mu_k(\sigma^2)$ is the k-th order central moment of a Gaussian law with variance $\sigma^2$.  Note that to compute the above integral, we performed a change of variable   $z-z^{M*}=\sqrt{\e}\,y$. Therefore each term $z-z^{*}$ can be considered as of order $\sqrt{\e}$ in the integration. This is why, to obtain a first order approximation of the moments in terms of $\e$, it is enough to have a fourth order approximation of $u(z)$, a second order approximation of $v_i(z)$ and a zero order approximation of $w_i(z)$,  in terms of $z$ around $z^{*}$.
The above computation leads in particular to the following approximations of the population size, the mean, the variance and the skewness of the population's distribution:
$$
\begin{cases}
N_{\e,i} =\int n_{\e,i} (z)dz=N_i^{M*}(1+\e (F_i+\f{E_i+0.5D_i^2}{A}+\f{3(C+BD_i)}{A^2}+\f{7.5B^2}{A^3}))+O(\e^2),
\\
\mu_{\e,i} =\f{1}{N_{\e,i} }\int z  n_{\e,i} (z) dz=z^{M*}+\e(\f{3B}{A^2}+\f{D_i}{A} )+O(\e^2),
\\
\sigma_{\e,i}^{  2}=\f{1}{N_{\e,i} }\int (z-\mu_{\e,i} )^2  n_{\e,i} (z) dz =\f{\e}{A} +O(\e^2),
\\
s_{\e,i} =\f{1}{\sigma_{\e,i}^{  3}N_{\e,i} }\int (z-\mu_{\e,i} )^3  n_{\e,i} (z) dz =\f{6B}{A^{\f32}}\sqrt{\e}+O(\e^{\f32}).
\end{cases}
$$

\section{Identification of the  ESS  (the proofs of  Theorem \ref{th:ESS} and Proposition \ref{prop:invade} )}

\label{sec:pr-ESS}

In this section, we prove Theorem \ref{th:ESS}, Corollary \ref{cor:deg} and Proposition \ref{prop:invade}. We first provide a description of the ESS in Subsection \ref{sec:ESS-des}. Next, we prove Theorem \ref{th:ESS}-(i) in Subsection \ref{sec:dim}. In Subsection  \ref{sec:mono} we prove  Theorem \ref{th:ESS}-(ii) and Corollary \ref{cor:deg}.  Finally in Subsection \ref{sec:inv-con} we prove Proposition \ref{prop:invade}.  \\

\subsection{The description of the ESS}
\label{sec:ESS-des}
\noindent
We first rewrite the conditions for ESS in terms of the following variables:
\beq
\label{mu}
\mu_i (N_i)= \f{\kappa_i N_i + m_i -r_i }{g_i}, \qquad i=1,2,
\eeq
where $\mu_i$ is an indicator of the size of the population in patch $i$. 
In several parts of this paper, we will express the effective fitness as a function of $\mu_i$ instead of $N_i$:
$$
W_\mu(z,\mu_1(N_1),\mu_2(N_2)) = W(z,N_1,N_2),
$$
hence, the effective fitness in terms of $\mu_i$  is given by
$$
\begin{array}{rl}
W_\mu(z,\mu_1,\mu_2)&=\f 12\left[ -g_1 (\mu_1+(z+\theta)^2)-g_2(\mu_2+(z-\theta)^2 ) \right.\\
&\left.+ \sqrt{ \big( g_1 (\mu_1+(z+\theta)^2) - g_2 (\mu_2+(z-\theta)^2) \big)^2+4m_1m_2 } \right].
\end{array}
$$
From the definition of ESS, we deduce that at the demographic equilibrium of an ESS, where the indicators of population size in patches $1$ and $2$ are given by $(\mu_1^*,\mu_2^*)$, we have
$$
W_\mu(z,\mu_1^*,\mu_2^*)\leq 0,\qquad \text{for $z\in \R$},
$$
with the equality attained at one or two points corresponding to the monomorphic or dimorphic ESS.  We then notice that the above inequality is equivalent with
$$
\begin{cases}
g_1(\mu_1^*+ (z+\theta)^2)+ g_2 (\mu_2^*+(z-\theta)^2) \geq 0,\\
f(z;\mu_1^*,\mu_2^*):=(\mu_1^*+(z+\theta)^2) (\mu_2^*+(z-\theta)^2) \geq \f {m_1m_2} {g_1g_2}.
\end{cases}
$$
This implies that at the ESS, $\mu_i^*>0$ and 
\beq
\label{min-pb}
\min_x \; (\mu_1^*+(z+ \theta)^2) (\mu_2^* +(z-\theta)^2) =\f{m_1 m_2}{g_1g_2}.
\eeq
Note that the above function is a fourth order polynomial and hence has one or two minium points, which here will correspond to the monomorphic or dimorphic ESS.  Conditions for the demographic equilibria will help us determine $(\mu_1^*,\mu_2^*)$:\\

\noindent
(i) If the minimum in \fer{min-pb} is attained at the point $z^{M*}$, for $z^{M*}$ to be an ESS the following condition must be satisfied:
$$
\left(
\begin{array}{cc}
-g_1\Big((z^{M*}+ \theta)^2+\mu_1^{*} \Big) & m_2 \\
m_1 & -g_2\Big((z^{M*}- \theta)^2+\mu_2^{*} \Big)
\end{array}
\right)
\left(
\begin{array}{c}
N_1^{M*}
\\
N_2^{M*}
\end{array}
\right)
=0,
$$
with
$$
N_i^{M*}>0,\qquad \mu_i^*=\mu_i (N_i^{M*})= \f{\kappa_i N_i^{M*} + m_i -r_i }{g_i}, \qquad i=1,2.
$$


\noindent 
(ii) If the minimum in \fer{min-pb} is attained at two points $z_{\mathrm{I}}^{D*}$ and $z_{\mathrm{II}}^{D*}$, for $(z_{\mathrm{I}}^{D*},z_{\mathrm{II}}^{D*})$ to be an ESS, there must exist $\nu_{k,i}>0$, for $i=1,2$ and $k=\mathrm{I},\mathrm{II}$, such that,  
\beq
\label{dem-dim-1}
\left(
\begin{array}{cc}
-g_1\Big((z_k^{D*}+ \theta)^2+\mu_1^{*} \Big) & m_2 \\
m_1 & -g_2\Big((z_k^{D*}- \theta)^2+\mu_2^{*} \Big)
\end{array}
\right)
\left(
\begin{array}{c}
\nu_{k,1}
\\
\nu_{k,2}
\end{array}
\right)
=0,
\qquad k=\mathrm{I},\mathrm{II},
\eeq
\beq
\label{dem-dim-2}
\nu_{\mathrm{I},1}+\nu_{\mathrm{II},,1}=N_1^{D*}, \quad \nu_{\mathrm{I},,2}+\nu_{\mathrm{II},,2}=N_2^{D*}, \quad \mu_i^*=\mu_i(N_i^{D*}) \; \text{ for $i=1,2$}.
\eeq

%
%

\subsection{The dimorphic ESS}
\label{sec:dim}

To identify the dimorphic ESS we first give the following lemma
\begin{lemma}
\label{lem:min}
If $f(z ; \mu_1,\mu_2)$ has two global minimum points $z_{\mathrm{I}}$ and $z_{\mathrm{II}}$, then $\mu_1=\mu_2$ and $z_{\mathrm{I}}=-z_{\mathrm{II}}$. 
\end{lemma}

\proof
Let's suppose that $f(z ; \mu_1,\mu_2)$ has two global minimum points $z_{\mathrm{I}}$ and $z_{\mathrm{II}}$ and $\mu_2<\mu_1$. The case with $\mu_1<\mu_2$ can be treated following similar arguments. \\

\noindent 
Since $z_{\mathrm{I}}$ and $z_{\mathrm{II}}$ are minimum points we have
$$
(\mu_1+(z_k+\theta)^2) (\mu_2 +(z_k-\theta)^2) \leq (\mu_1+(-z_k+\theta)^2) (\mu_2 +(-z_k-\theta)^2), \qquad k=\mathrm{I},\mathrm{II}.
$$
It follows that
$$
0 \leq 4z_k\theta (\mu_1-\mu_2), \qquad k=\mathrm{I},\mathrm{II},
$$
and hence 
$$
0\leq z_k, \qquad k=\mathrm{I},\mathrm{II}.
$$
This implies in particular that all the roots of $f'(z,\mu_1,\mu_2)$ are positive. However, this is not possible since
$$
f'(z,\mu_1,\mu_2) = 4z^3+2(\mu_1+\mu_2-2\theta^2)z+2\theta(\mu_2-\mu_1).
$$
The fact that there is no second order term in the above expression implies that the sum of the roots is zero and hence the roots change sign. This is a contradiction with the previous arguments. We hence deduce that $\mu_1=\mu_2$.
\qed

\noindent
The above lemma indicates that at a dimorphic ESS one should have $\mu_1^*=\mu_2^*=\mu^*$. Hence to find a dimorphic ESS we look for $(\mu^*,z_{\mathrm{I}}^*,z_{\mathrm{II}}^*)$ such that
\beq
\label{min-dim}
f(z_k^*,\mu^*,\mu^*)=\min f(z;\mu^*,\mu^*)=\f {m_1m_2}{g_1g_2}, \qquad k=\mathrm{I},\mathrm{II}.
\eeq
To identify the minimum points of $f$ we differentiate $f$ with respect to $z$ and find
$$
f'(z,\mu^*,\mu^*) = 4z^3+4(\mu^* -\theta^2)z.
$$
For $f$ to have two minimum points, $f'$ must have three roots and hence one should have
\beq
\label{cond-mu}
\mu^* < \theta^2.
\eeq
Then, the minimum points are given by
$$
z_{\mathrm{I}}^*=-\sqrt{\theta^2-\mu^*}, \qquad z_{\mathrm{II}}^*=\sqrt{\theta^2-\mu^*}.
$$
Then replacing the above values in \fer{min-dim} we obtain
$$
\mu^*= \f{m_1m_2}{4\theta^2g_1g_2}.
$$
Note that combining the above line with condition \fer{cond-mu} we obtain  \fer{as1:dim}. \\

\noindent
Up until now, we have proven that if a dimorphic ESS exists  \fer{as1:dim} is verified and the dimorphic ESS is given by $(z_{\mathrm{I}}^{D*},z_{\mathrm{II}}^{D*})=(-\sqrt{\theta^2-\mu^*},\sqrt{\theta^2-\mu^*})$. However, for this point to be an ESS, as explained in the previous subsection, there must exist  $\nu_{k,i}>0$, for $i=1,2$ and $k=\mathrm{I},\mathrm{II}$ such that \fer{dem-dim-1}--\fer{dem-dim-2} are satisfied. 
Replacing $z_k^{D*}$ by their values and solving \fer{dem-dim-1}--\fer{dem-dim-2}, we obtain that   $\nu_{k,i}$, for $i=1,2$ and $k=\mathrm{I},\mathrm{II}$, are identified in a unique way by \fer{nuij}. One can verify  by simple computations that the weights  $\nu_{k,i}$  are positive if and only if conditions \fer{as2:dim}--\fer{as3:dim} are satisfied. 
As a conclusion, we obtain that a dimorphic ESS exists if and only if the conditions \fer{as1:dim}--\fer{as3:dim} are satisfied. Moreover, when it exists, such dimorphic ESS is unique.

\subsection{The monomorphic ESS}
\label{sec:mono}

In this subsection we prove Theorem \ref{th:ESS}-(ii) and Corollary \ref{cor:deg}. To this end, we assume thanks to \fer{as:r-m} and without loss of generality that $r_1-m_1>0$ and then we consider two cases:\\

\noindent 
(i) We first suppose that condition \fer{as1:dim} does  not hold.
We  then introduce the following functions:
$$
\begin{cases}
F=(F_1,F_2):(0,+\infty) \to (0,+\infty) \times [-\theta, \theta]
\\
\mu_2 \mapsto (\mu_1, \overline z),
\end{cases}
\qquad
\begin{cases}
G :  (0,+\infty) \times  [-\theta, \theta] \to \R
\\
(\mu_1, \overline z) \mapsto \overline \mu_2,
\end{cases}
$$
where $\mu_1$ and $\overline z$ are chosen such that
$$
f(\overline z,\mu_1,\mu_2) = \min f(z ; \mu_1, \mu_2) =\f{m_1m_2}{g_1g_2},
$$
and $\overline \mu_2$ is  given by
\beq
\label{mu2bar}
\overline \mu_2 =\f{1}{g_2} \left[ \f{\kappa_2g_1}{m_2} ((\overline z + \theta)^2+\mu_1) \Big( \f{g_1\mu_1+r_1-m_1}{\kappa_1} \Big)+m_2-r_2 \right].
\eeq
We claim the following lemma which we will prove at the end of this paragraph.
\begin{lemma}
\label{lem:FG}
If \fer{as1:dim} does not hold, then the functions $F$ and $G$ are well-defined. Moreover, $F_1$ and $F_2$ are   decreasing with respect to $\mu_2$ and $G$ is increasing with respect to $\mu_1$ and $\overline z$.
\end{lemma}
Following the arguments in Section \ref{sec:ESS-des}, one can  verify that a trait $z^{*}$ is a monomorphic ESS with a demographic equilibrium $(\mu_1^*,\mu_2^*)$ if and only if
$F(\mu_2^*)=(\mu_1^*,z^*)$ and $G \circ F (\mu_2^*)=\mu_2^*$. Therefore, identifying monomorphic evolutionary stable strategies is equivalent with finding the fixed points of $G \circ F$. \\
\noindent 
In the one hand, from Lemma \ref{lem:FG} we deduce that $G \circ F$ is a decreasing function . In the other hand, one can verify that, as $\mu_2\to 0$, $G \circ F(\mu_2) \to +\infty$. In particular $G \circ F(\mu_2) > \mu_2$ for $\mu_2$ small enough.  It follows that there exists a unique $\mu_2^*$ such that $G \circ F(\mu_2^*)=\mu_2^*$. We deduce that there exists a unique ESS which is given by $z^{M*}=F_2(\mu_2^*)$. Moreover, $(F_1(\mu_2^*),\mu_2^*)$ corresponds to its demographic equilibrium.

\noindent
Note that for such ESS to make sense, one should also have $N_i^{M*}(\mu_i^*)>0$. This is always true for such fixed point. Note indeed that,  since 
$\mu_1^*=F_1(\mu_2^*)\in (0,\infty)$ and  $r_1-m_1>0$ we deduce that $N_1^{M*}>0$. Moreover, the positivity of $N_2^{M*}$ follows from $N_2^{M*}=(g_2\mu_2^*+r_2-m_2)/\kappa_2$, \fer{mu2bar} and the positivity of $r_1-m_1$ and $\mu_1^*$.
\\

{\bf Proof of Lemma \ref{lem:FG}.}
The fact that $G : (0,+\infty) \times  [-\theta, \theta] \to \R$ (and respectively $F_1=(0,+\infty) \to (0,\infty)$) is well-defined and increasing (respectively decreasing)   is immediate. 
We only show that $F_2$ is well-defined and decreasing.
To this end, we notice that since $f$ is a fourth order polynomial, it admits one or two minimum points. However, from the arguments in Subsection \ref{sec:dim} we know that the only possibility to have two global minima is that \fer{as1:dim} holds and $\mu_2=\mu^*=\f{m_1m_2}{4\theta^2g_1g_2}$. Since we assume that \fer{as1:dim} does not hold, $f$ always admits a unique minimum point in $\R$. This minimum point is indeed attained in $[-\theta,\theta]$ since for all $z<-\theta$, $f(z ; \mu_1,\mu_2) >f(-\theta; \mu_1,\mu_2)$ and
 for all $z>\theta$, $f(z ; \mu_1,\mu_2) > f( \theta; \mu_1,\mu_2)$. Hence $\overline z$ is defined in a unique way in $[-\theta, \theta]$.\\

\noindent
Finally, it remains to prove that $F_2 :(0,\infty) \to [-\theta,\theta]$ is a decreasing function. To this end, let's suppose that $\widetilde \mu_2>   \mu_2$. Therefore, $F_1(\widetilde \mu_2)=\widetilde \mu_1 < F_1(\mu_2)=\mu_1$. We want to prove that $F_2(\widetilde \mu_2)=\widetilde z <F_2(\mu_2)=\overline z$. To this end, we write
$$
\begin{array}{rl}
f(z ; \widetilde \mu_1,\widetilde \mu_2) &= f(z ; \mu_1,\mu_2) +(\widetilde  \mu_1- \mu_1) (z-\theta)^2  + (\widetilde  \mu_2- \mu_2) (z+\theta)^2 +\widetilde  \mu_1\widetilde  \mu_2- \mu_1\mu_2\\
& =f(z;\mu_1,\mu_2) +h(z; \mu_1,\mu_2, \widetilde \mu_1,\widetilde \mu_2),
\end{array}
$$
where $h$ is increasing with respect to $z$. Since $f(z,\mu_1,\mu_2)$ attains its minimum at $\overline z$ and $f(z, \widetilde \mu_1,\widetilde \mu_2)$ attains its minimum at $\widetilde  z$ we find that
$$
f(\overline z ; \mu_1,\mu_2) < f(\widetilde z ; \mu_1,\mu_2),
$$
$$
f(\widetilde z ; \mu_1,\mu_2) + h(\widetilde z ; \mu_1,\mu_2, \widetilde \mu_1,\widetilde \mu_2) < f(\overline z ; \mu_1,\mu_2) + h(\overline z ; \mu_1,\mu_2, \widetilde \mu_1,\widetilde \mu_2).
$$
Combining the above inequalities, we obtain that
$$
h(\widetilde z ; \mu_1,\mu_2, \widetilde \mu_1,\widetilde \mu_2) <  h(\overline z ; \mu_1,\mu_2, \widetilde \mu_1,\widetilde \mu_2).
$$
and since $h$ is an increasing function, we conclude that
$
\widetilde z < \overline z.
$

\qed

\bigskip

\noindent 
(ii) We next suppose that \fer{as1:dim}   holds. Consequently, $F$ is not well-defiled at $\mu_2=\mu^*=\f{m_1m_2}{4\theta^2g_1g_2}$ since $F_1(\mu^*)=\mu^*$ and $\max_z f(z;\mu^*,\mu^*)$ is attained at two points $\pm z^{D*}$. Therefore, we only can define $F$ in $(0,\infty)\setminus \{\mu^*\}$:
$$
\begin{cases}
\widetilde F=(\widetilde F_1,\widetilde F_2):(0,+\infty)\setminus \{\mu^*\} \to (0,+\infty) \times [-\theta, \theta]
\\
\mu_2 \mapsto (\mu_1, \overline z),
\end{cases}
\qquad
\begin{cases}
G :  (0,+\infty) \times  [-\theta, \theta] \to \R
\\
(\mu_1, \overline z) \mapsto \overline \mu_2,
\end{cases}
$$
where $\mu_1$, $\overline z$ and $\overline m_2$ are chosen as above. Following similar arguments as in the proof of Lemma \ref{lem:FG} we obtain
\begin{lemma}
\label{lem:FG2}
Under condition \fer{as1:dim} the functions $\widetilde F$ and $G$ are well-defined. Moreover, $\widetilde F_1$ and $\widetilde F_2$ are   decreasing with respect to $\mu_2$ in the intervals $(0,\mu^*)$ and $(\mu^*,+\infty)$ and $G$ is increasing with respect to $\mu_1$ and $\overline z$.
\end{lemma}
As above, identifying monomorphic evolutionary stable strategies is equivalent with finding the fixed points of $G \circ \widetilde F$, which is a decreasing function  in the intervals $(0,\mu^*)$ and $(\mu^*,+\infty)$  thanks to the  lemma \ref{lem:FG2}. We then compute
$$
\widetilde F(\mu^{*-})=(\mu^*, z^{D*}),\qquad F(\mu^{*+})=(\mu^*,-z^{D*}),
$$
$$
G \circ \widetilde F(\mu^{*-}) =\f{1}{g_2} \left[ \f{\kappa_2}{m_2} (g_1( z^{D*} + \theta)^2)+\mu^*) \Big( \f{g_1\mu^* +r_1-m_1}{\kappa_1} \Big)+m_2-r_2 \right],
$$
$$
G \circ \widetilde F(\mu^{*+}) =\f{1}{g_2} \left[ \f{\kappa_2}{m_2} (g_1( -z^{D*} + \theta)^2)+\mu^*) \Big( \f{g_1\mu^* +r_1-m_1}{\kappa_1} \Big)+m_2-r_2 \right],
$$
where $\mu^{*+}$and $\mu^{*-}$ correspond respectively to the limits from the right and from the left as $\mu \to \mu^*$.
One can easily verify that $G \circ \widetilde F(\mu^{*+}) < \mu^*$ if and only if \fer{as2:dim} holds, and similarly $G \circ \widetilde F(\mu^{*-}) > \mu^*$ if and only if 
\fer{as:eq8}, or equivalently
 \fer{as3:dim},  holds. We hence deduce, from the latter property and the fact that $G \circ \widetilde F$ is decreasing in the intervals $(0,\mu^*)$ and $(\mu^*,+\infty)$, that:
\begin{enumerate}
\item If \fer{as2:dim} and \fer{as3:dim} hold there is no monomorphic ESS. Note that, under these conditions there exists a unique dimorphic ESS.

\item   if  \fer{as2:dim} holds and  the r.h.s. of \fer{as3:dim} is negative, then there exists a unique monomorphic ESS in $\mu_2^{M*}\in (0,\mu^*)$, $\mu_1^{M*}\in (\mu^*,\infty)$ and $z^{M*}\in (z^{D*},\theta)$.

\item if  \fer{as3:dim} holds and   the r.h.s. of  \fer{as2:dim} is negative, then there exists a unique monomorphic ESS with $\mu_2^{M*}\in (\mu^*,\infty)$, $\mu_1^{M*}\in (0,\mu^*)$ and $z^{M*}\in (-\theta,-z^{D*})$.

\item   if  \fer{as2:dim} holds and   \fer{as3:dim} is an equality, then there exists a unique monomorphic ESS which is given by $\{z^{D*}\}$ and $\mu_1^*=\mu_2^*=\mu^*$. 

\item   if  \fer{as3:dim} holds and   \fer{as2:dim} is an equality, then there exists a unique monomorphic ESS which is given by $\{-z^{D*}\}$ and $\mu_1^*=\mu_2^*=\mu^*$. 

\item Finally, from the fact that \fer{as2:dim} and \fer{as3:dim} are respectively equivalent to \fer{as:eq7} and \fer{as:eq8} we deduce that at least one of conditions \fer{as2:dim} and \fer{as3:dim} always holds. Therefore, all the possible cases have been considered.
\end{enumerate}

\noindent
Note that, following similar arguments to the previous case, the total population sizes $N_i^{M*}(\mu_i^*)$, for $i=1,2$, corresponding to the unique fixed point, are positive and hence the obtained monomorphic ESS is indeed valid.  This concludes the proof of Theorem \ref{th:ESS}. It remains to prove Corollary \ref{cor:deg}:\\

{\bf Proof of Corollary \ref{cor:deg} } We first notice from the arguments above that $W(z,N_1^*,N_2^*)=W_\mu(z,\mu_1^*,\mu_2^*)$ has at most two global maximum points. Therefore, for \fer{Wneg} not to hold, the unique ESS should be monomorphic while $W_\mu(z,\mu_1^*,\mu_2^*)$ has two maximum points.
However, from the arguments in Section \ref{sec:dim} we know that if $W_\mu(z,\mu_1^*,\mu_2^*)$ has two maximum points, then \fer{as1:dim}  holds, $\mu_1^*=\mu_2^*=\mu^*$ and the maximum points are given by $\{\pm z^{D*}\}$. Finally, from the results in the above paragraph, we know that the only possibility to have a monomorphic ESS in this case, is that either \fer{as2:dim} or \fer{as3:dim} is an equality, which is in contradiction with \fer{non-deg}.
\qed

\subsection{The interpretation of conditions \fer{as2:dim} and \fer{as3:dim} }
\label{sec:inv-con}

In this subsection we prove Proposition \ref{prop:invade}. We only prove the first claim. The second claim can be derived following similar arguments.\\

\noindent
We denote by $(\mu_1^{\rm eq},\mu_2^{\rm  eq})$ the demographic equilibrium  of a monomorphic population of trait $-z^{D*}$ and we first claim the following lemma.
\begin{lemma}
There exists a unique demographic equilibrium $n_i=N_i \delta (z+z^{D*})$ corresponding to the the set $\Omega=\{-z^{D*}\}$.
\end{lemma}

\proof
We introduce two functions $K$ and $H$ which are respectively close to $F_1$ and $G$ introduced above:
$$
\begin{cases}
K:(-(z^{D*}+\theta)^2,+\infty) \to \R 
\\
\mu_2 \mapsto \mu_1,
\end{cases}
\qquad
\begin{cases}
H :  \R  \to \R
\\
\mu_1  \mapsto \overline \mu_2,
\end{cases}
$$
where $\mu_1$ is chosen such that
$$
f(-z^{D*} ; \mu_1,\mu_2)=\f{m_1m_2}{g_1g_2},
$$
and $\overline \mu_2$ is  given by
$$
\overline \mu_2 =\f{1}{g_2} \left[ \f{\kappa_2 g_1}{m_2} ((z^{D*} - \theta)^2+\mu_1) \Big( \f{g_1\mu_1+r_1-m_1}{\kappa_1} \Big)+m_2-r_2 \right].
$$
Then  the demographic equilibrium  $(\mu_1^{\rm eq}, \mu_2^{\rm  eq})$ of a monomorphic resident population of type $-z^{D*}$ corresponds to a fixed point of $H \circ K$: 
$$
H \circ K(\mu_2^{\rm  eq}) =\mu_2^{\rm  eq}, \qquad 
K(\mu_2^{\rm  eq})= \mu_1^{\rm  eq}.
$$
Note also that, for such equilibrium to make sense, one should have $0 \leq N_i(\mu_i^{\rm eq}) $ or equivalently
$$
\f{m_i-r_i}{g_i} \leq \mu_i^{\rm eq} .
$$
Moreover, since $W_\mu(-z^{D*}, \mu_1^{\rm eq}, \mu_2^{\rm eq})=0$, we have the additional condition 
$$
0<\mu_1^{\rm eq}+(z^{D*}-\theta)^2, \qquad 0<\mu_2^{\rm eq}+(z^{D*}+\theta)^2.
$$
Reciprocally, a pair $(\mu_1,\mu_2)$ which satisfies the above conditions corresponds to a demographic equilibrium.\\

\noindent
We next notice, on the one hand, that $K$ is a decreasing function, and hence, in view of the above conditions,  a fixed point $(\mu_1^{\rm eq}, \mu_2^{\rm eq})$ of $H \circ K$, is a demographic equilibrium if and only if $\mu_2^{\rm eq}\in \big(   -(-z^{D*}+\theta)^2, \widetilde \mu_2)$, with $\widetilde \mu_2=K^{-1}(\max ( \f{m_1-r_1}{g_1},-( z^{D*}-\theta)^2 )$. On the other hand, 
 $H$, restricted to $\big(\max ( \f{m_1-r_1}{g_1},-( z^{D*}-\theta)^2 ) ,+\infty\big)$, is an increasing function. Therefore $H \circ K$, restricted to the set $\big(   -(z^{D*}+\theta)^2, \widetilde \mu_2)$,  is decreasing.  We deduce that a demographic equilibrium, if it exists, is unique.\\
 
 \noindent
We then note that, as $\mu_2 \to -(z^{D*}+\theta)^{2+}$, $H \circ K(\mu_2) \to +\infty$. In particular, for $\mu_2$ close to $ -(z^{D*}+\theta)^{2}$, $H \circ K(\mu_2)>\mu_2$. Furthermore, $H \circ K( \widetilde \mu_2)=\f{m_2-r_2}{g_2}<0$. Note also that, $K(\widetilde \mu_2)<0$ and $K(\mu^*)=\mu^*>0$ and hence $0<\mu^*<\widetilde \mu_2$, which implies that $H \circ K( \widetilde \mu_2)<\widetilde \mu_2$. We deduce from the intermediate value theorem that, $H\circ K: \big(   -(z^{D*}+\theta)^2, \widetilde \mu_2)\to \R$ has a unique fixed point $(\mu_1^{\rm eq}, \mu_2^{\rm eq})$ and hence there exists a unique demographic equilibrium.
\qed

\medskip


\noindent 
We next observe that, since  $W_{\mu}(-z^{D*},\mu_1^{\rm eq},\mu_2^{\rm eq})=0$,  we have $W_{\mu}(z^{D*},\mu_1^{\rm eq},\mu_2^{\rm eq})>0$ if and only if $\mu_2^{\rm eq}<  \mu_1^{\rm eq}$. Moreover, since $W_\mu(-z^{D*},\mu^*,\mu^*)=0$, this is equivalent with $\mu_2^{\rm eq}< \mu^*< \mu_1^{\rm eq}$.  
\\

\noindent
We are now ready to conclude. Let's first suppose that \fer{as2:dim} holds which implies that $H \circ K(\mu^*)<\mu^*$. Then, thanks to the fact that $\mu^*<\widetilde \mu_2$ and from the monotonicity of $K$ and $H \circ K$ we deduce that the unique fixed point, $\mu_2^{\rm eq}$, of $H \circ K$  satisfies 
$$
\mu_2^{\rm eq}<\mu^* <K(\mu_2^{\rm eq})=:\mu_1^{\rm eq}.
$$
This implies that $W_{\mu}(z^{D*},\mu_1^{\rm eq},\mu_2^{\rm eq})>0$ or equivalently, a mutant trait $z^{D*}$ can invade a resident population of trait $-z^{D*}$ at its demographic equilibrium. \\

\noindent
Let's now suppose that $W_{\mu}(z^{D*},\mu_1^{\rm eq},\mu_2^{\rm eq})>0$ and hence $\mu_2^{\rm eq}< \mu^*< \mu_1^{\rm eq}$. We then deduce from $H \circ K(\mu_2^{\rm eq})=\mu_2^{\rm eq}$ and that the monotonicity of $H \circ K$ that $H \circ K(\mu^*)<\mu^*$. This implies \fer{as2:dim}.

\section{The proof of Theorem \ref{thm:main}}
\label{sec:pr-HJ}

In this section, we prove Theorem \ref{thm:main}. To this end, we first provide a convergence result along subsequences in Subsection \ref{sec:HJ}. We next conclude using a uniqueness argument in Subsection \ref{sec:HJ-unique}.

\subsection{Convergence to the Hamilton-Jacobi equation with constraint }
\label{sec:HJ}

In this section, we prove that as $\e\to 0$, both sequences $(u_{\e,i})_\e$, for $i=1,2$, converge along subsequences and locally uniformly to a  function $u\in \mathrm{C}(\R)$,  such that $u$ is a viscosity solution to the following equation
\begin{equation}
\label{HJ-u}
\left\{ \begin{array}{ll}-| u'(z)|^2= W(z,N_1,N_2),&\quad \text{in $\R$},\\
\max_{z\in \R}u(z)=0,\end{array}\right.
\end{equation}
\beq
\label{supp-n}
{\rm supp}\, n_1= {\rm supp}\, n_2\subset \{z\, |\, u(z)=0 \} \subset \{z\, |\, W(z,N_1,N_2)=0\},
\eeq
where $(n_1,n_2)$ (respectively $(N_1,N_2)$) is a limit, along subsequences, of $(n_{\e,1},n_{\e,2})$ (respectively $(N_{\e,1},N_{\e,2})$) as $\e$ vanishes. Moreover, 
$$
N_i=\int_\R n_i(z)dz.
$$
Note that this is indeed the claim of Theorem \ref{thm:main}, except that we don't  know yet if $(n_1,n_2)=(n_1^*,n_2^*)$.\\
To this end, we first claim the following

\begin{prop}
\label{prop:Ne} Assume \fer{as:r-m}--\fer{as:m}. \\
(i) For all $\e>0$, we have
\beq
\label{N-bound}
N_{\e,1}+N_{\e,2} \leq N_M=2 \max (r_1,r_2).
\eeq
In particular, for $i=1,2$, $(n_{\e,i})_\e$ converge along subsequences  and weakly in the sense of measures to $n_i$ and $N_{\e,i}$ converges along subsequences to $N_i$.\\
(ii) For any compact set $K\subset \R$, there exists a constant $C_M=C_M(K)$ such that, for all $\e\leq 1$,
\beq
\label{Harnack}
  n_{\e,i}(x) \leq C_M n_{\e,j}(y), \qquad \text{for}\quad i,j\in \{1,2\},  |x-y|\leq \e.
\eeq
(iii)  For all $\eta>0$ there exists a constant $R$ large enough such that 
\beq
\label{n-da}
\int_{|z|>R} n_{\e,i}(z)dz <\eta, \quad \text{for }i=1,2.
\eeq
Consequently $N_i=\int_\R n_i(z) dz$. 

\end{prop}
We postpone the proof of this proposition to the end of this paragraph and we pursue giving the scheme of the proof of Theorem \ref{thm:main}. The next step, is to introduce functions $(l_{\e,1},l_{\e,2})$ as below
$$
l_{\e,i}:=\alpha_\e n_{\e,i}, \quad \text{for }i=1,2,
$$
with $\alpha_\e$ chosen such that
\beq
\label{int-l}
\int_\R \big( l_{\e,1}(z)+l_{\e,2}(z) \big)dz =1.
\eeq
Moreover, we define
$$
v_{\e,i}:= \e \log (l_{\e,i}),\quad \text{for }i=1,2.
$$
We next prove the following
\begin{prop}
\label{prop:ve}
Assume \fer{as:r-m}--\fer{as:m}. \\
(i) For $i=1,\,2$ and all $\e\leq \e_0$, the families $(v_{\e,i})_\e$ are locally uniformly bounded and locally uniformly Lipschitz.\\
(ii) As $\e\to 0$, both families $(v_{\e,i})_\e$, for $i=1,2$, converge along subsequences and locally uniformly in $\R$ to a continuous function $v\in \mathrm{C}(\R)$ and  $(N_{\e,i})_\e$, for $i=1,2$, converge along  subsequences  to $N_i$,  such that $v$ is a viscosity solution to the following equation
\begin{equation}
\label{HJ-v}
\left\{ \begin{array}{ll}-| v'(z)|^2= W(z,N_1,N_2),&\quad \text{in $\R$},\\
\max_{z\in \R}v(z)=0.\end{array}\right.
\end{equation}
(iii)  We have 
\beq
\label{W-neg}
W(z,N_1,N_2) \leq 0.
\eeq
Consequently, there exists $\delta >0$ such that
\beq
\label{bd-bel-N}
N_i \geq \delta, \quad \text{for }i=1,2.
\eeq
\end{prop}
The proof of this proposition is given at the end of this subsection. Note that \fer{bd-bel-N} implies that, for $\e$ small enough, $N_{\e,i}\geq \f \delta 2$. This together with \fer{N-bound} imply that, for $\e\leq \e_1$ with $\e_1$ small enough,
$$
\f{1}{2\max(r_1,r_2)} \leq \alpha_\e  \leq \f 1 \delta,
$$
and consequently
$$
 v_{\e,i}+\e \log (\delta)    \leq u_{\e,i} \leq v_{\e,i}+\e \log (2\max(r_1,r_2)).
$$
We then conclude from the above inequality together with  Proposition \ref{prop:ve}--(ii) that $(u_{\e,i})_\e$, for $i=1,2$, converge along subsequences and locally uniformly to a  function $u\in \mathrm{C}(\R)$ which is a viscosity solution of \fer{HJ-u}. \\
To prove \fer{supp-n}  we use the following lemma:
\begin{lemma}
\label{lem:sconv}
The function $v$ is semiconvex.
\end{lemma}
Then \fer{supp-n} is immediate from the WKB ansatz \fer{WKB} and the fact that $v$ is differentiable at its maximum points (since  it is a semiconvex function). Finally, lemma \ref{lem:sconv} can be proved following similar arguments as in  \cite{SM:12}--Theorem 1.2,  but  using  cut-off functions to treat the unbounded case as in the proof of Proposition \ref{prop:ve}-(i).\\

{\bf Proof of Proposition \ref{prop:Ne}.}
(i) We first prove \fer{N-bound}. To this end, we integrate the equations in \fer{main} with respect to $z$ to obtain
$$
\int_\R n_{\e,i}(z)(r_i-m_i-g_i(z-\theta_i)^2-N_{\e,i}) dz +m_j N_{\e,j}=0,\quad i=1,2, \; j=2,1.
$$
Adding the two equations above, it follows that
$$
N_{\e,1}^2+N_{\e,2}^2 \leq r_1 N_{\e,1}+r_2 N_{\e,2},
$$ 
and hence \fer{N-bound}.

(ii) We define 
$$
K_\e=\Big\{\f x \e \,  |\, x\in K \Big\},  \qquad  \widetilde n_{\e,i}(y)=n_{\e,i}(\e y),\qquad\text{for $i=1,\,2$}.
$$
From \fer{main} we have, for $z\in \R$,
\begin{equation}
\label{tilden}
\displaystyle\left\{\begin{array}{rll}
-  \widetilde n_{\e,1}''(z)=\widetilde n_{\e,1}(z) R_1(\e z,N_{\e,1})+m_2 \widetilde n_{\e,2}(z)-m_1\widetilde n_{\e,1}(z),\\
-  \widetilde n_{\e,2}''(z)=\widetilde n_{\e,2}(z) R_2(\e z,N_{\e,2})+m_1\widetilde n_{\e,1}(z)-m_2\widetilde n_{\e,2}(z).
\end{array}\right.
\end{equation}
Moreover, from \fer{Ri} and \fer{N-bound} we obtain that there exists a constant $C=C(K)$ such that
$$
-C \leq R_i(\e z,{N_{\e,i}})\leq C, \qquad \text{for all } z\in K_{\e}.
$$
Therefore the coefficients of the linear elliptic system \fer{tilden} are bounded uniformly in $K_\e$. It follows from the classical Harnack inequality (\cite{JB.MS:04}, Theorem 8.2) that there exists a constant $C_M=C_M(K)$ such that, for all $z_0\in K_\e$ such that $B_1(z_0)\subset K_\e$ and for $i,j=1,2$,
$$
\sup_{z\in B_1(z_0)}\widetilde n_\e^i(z)\leq C_M\,\inf_{z\in B_1(z_0)}\widetilde n_\e^j(z).
$$
Rewriting the latter in terms of $ n_\e^1$ and $ n_\e^2$ and replacing $(z,z_0)$ by $(\frac{z'}{\e},\f{z'_0}{\e})$ we obtain
$$
\sup_{z'\in B_\e(z'_0)}n_\e^i(z')\leq C_M\,\inf_{z' \in B_\e(z'_0)}n_\e^j(z'),\quad
$$
and hence \fer{Harnack}.

(iii) We integrate the equations in \fer{main} with respect to $z$ to obtain 
\beq
\label{int-n}
0\leq \int_\R n_{\e,i}(z) (r_i-g_i(z+\theta)^2 )dz +m_j N_{\e,j}(z)  .
\eeq
We choose a constant $R>0$ large enough such that for all $|z|>R$, we have
$$
r_i-g_i(z-\theta_i)^2  < - \f{N_M}{\eta} \max(r_1+m_2,r_2 +m_1),\qquad i=1,2.
$$
Splitting the integral term in the r. h. s. of \fer{int-n} into two parts we obtain
$$
0< r_i \int_{|z|\leq R} n_{\e,i}(z)dz -\f{N_M}{\eta} \max(r_1+m_2,r_2 +m_1) \int_{|z|>R} n_{\e,i}(z)dz +m_j N_{\e,j}.
$$
Next, using \fer{N-bound}, we obtain
$$
\f{N_M}{\eta} \max(r_1+m_2,r_2 +m_1) \int_{|z|>R} n_{\e,i}(z)dz  < (r_i+m_j)N_M,
$$
and hence \fer{n-da}.\\
\qed

{\bf Proof of Proposition \ref{prop:ve}.}
(i) We first prove that for all $a>0$ and any compact set $K$, there exists $\e_0$ such that for all $\e\leq \e_0$, we have
$$
v_{\e,i}(z) \leq a,\quad \text{for  } i=1,2,\quad z\in K.
$$
Note that, thanks to \fer{Harnack}, for any compact set $K$, there exists a constant $C_M=C_M(K)$ such that
\beq
\label{v-Har}
|v_{\e,i}(x)-v_{\e,j}(y)| \leq \e \log C_M, \quad \text{for  $|x-y|\leq \e$ and $i=1,2$.}
\eeq
We fix a compact set $K$. Let $z_0\in K$, $i\in \{1,2\}$ and $\e \leq \e_0=\f{a}{2\log C_M}$ be such that
$$
a < v_{\e,i}(z_0).
$$
Therefore,  for all $|y-z_0|\leq \e$, we find
$$
\f a 2 < a-\e \log C_M <v_{\e,i}(y).
$$
It follows that
$$
\e \exp(\f{a}{2\e}) \leq   \int_{|y-z_0|\leq \e} \exp( \f{v_{\e,i}(y)}{\e} )dy \leq \int_\R l_{\e,i}(y)dy.
$$
Note that the l. h. s. of the above inequality tends to $+\infty$ as $\e\to 0$, while  the r. h. s. is bounded by $1$, which is a contradiction. Such $z_0$ therefore does not exists and for all $z\in K$, $\e\leq \e_0$ and $i=1,2$, we find
$$
v_{\e,i}(z)\leq a.
$$

(ii) We next notice that, similarly to the proof of Theorem \ref{prop:Ne}-(iii), one can prove that, for all $\eta>0$ there exists a constant $R$ large enough such that 
\beq
\label{l-da}
\int_{|z|>R} l_{\e,i}(z)dz <\eta, \quad \text{for }i=1,2.
\eeq

(iii) Next, we prove that there exists $\e_0>0$ such that for all $\e\leq \e_0$, the families $(v_{\e,i})_\e$ are locally uniformly bounded from below. To this end, we first observe from \fer{int-l} and \fer{l-da} that, for $\eta\in (0,\f 14)$  there exists a constant $R_0>0$ such that  
$$
\int_{|z|\leq R_0}( l_{\e,1}(z) +l_{\e,2}(z) )dz > 1-2\eta>\f 1 2.
$$
Consequently, for $\e\leq \e_0$, with $\e_0$ small enough, there exists $z_0\in \R$ and $i\in \{1,2\}$  such that $|z_0|\leq R_0$ and $-1\leq v_{\e,i}(z_0) $. We deduce, thanks to \fer{v-Har}, that for any compact set $K=\overline B_R(0)$, with $R\geq R_0$,
$$
-1-2\log(C_M(K))R \leq v_{\e,i} (z),\qquad \text{for }i=1,2, \; \e\leq \e_0, \; z\in K .
$$

(iv) We prove that, for any compact set $K$, the families $(v_{\e,i})_\e$ are  uniformly Lipschitz in $K$. To this end, we first notice that $(v_{\e,i})_\e$ solves the following system:
\beq
\label{sys-ve}
-\e v_{\e,i}''= |v_{\e,i}'|^2+R_i(z,N_{\e,i}) + m_j \exp \Big( \f{v_{\e,j}-v_{\e,i}}{\e} \Big) -m_i, \quad i=1,2, \; j=2,1.
\eeq
We differentiate the above equation with respect to $z$ and multiply it by $v_{\e,i}'$ to obtain
$$
-\e v_{\e,i}' v_{\e,i}'''=2v_{\e,i}'^2 v_{\e,i}''+ \f{\p}{\p z}R_i(z,N_{\e,i}) v_{\e,i}' +m_j v_{\e,i}' \Big(  \f{v_{\e,j}'-v_{\e,i}'}{\e} \Big) \exp \Big( \f{v_{\e,j}-v_{\e,i}}{\e} \Big).
$$
We then define $p_{\e,i} := |v_{\e,i}'|^2$ and notice that
$$
p_{\e,i}'=2 v_{\e,i}' v_{\e,i}'',\qquad p_{\e,i}''=2 v_{\e,i}''^2+2v_{\e,i}'v_{\e,i}'''.
$$
Combining the above lines we obtain that
\beq
\label{eq-p}
-\f \e 2 p_{\e,i}'' +\e v_{\e,i}''^2=2 p_{\e,i}'v_{\e,i}'+ \f{\p}{\p z}R_i(z,N_{\e,i}) v_{\e,i}' +m_j v_{\e,i}' \Big(  \f{v_{\e,j}'-v_{\e,i}'}{\e} \Big) \exp \Big( \f{v_{\e,j}-v_{\e,i}}{\e} \Big).
\eeq
We then fix a point $\xi \in K$ and introduce a cut-off function $\vp\in C^{\infty}(\R)$ which satisfes
\beq
\label{vp}
\vp(\xi)=1,\quad 0\leq \vp \leq 1 \text{ in $\R$}, \quad \vp\equiv 0 \text{ in $B_1(\xi)^c$},\quad |\vp'| \leq C \vp^{\f 12}, \qquad |\vp''|\leq C.
\eeq
We then define $P_{\e,i}=p_{\e,i}\vp$ and notice that
$$
P_{\e,i}'=p_{\e,i}' \vp+p_{\e,i} \vp',\qquad P_{\e,i}''= p_{\e,i}'' \vp+2 p_{\e,i}'\vp'+p_{\e,i}\vp''.
$$
We then multiply \fer{eq-p} by $\vp$ to obtain
$$
\begin{array}{rl}
-\f \e 2 P_{\e,i}'' +\e \vp v_{\e,i}''^2&=2 P_{\e,i}'v_{\e,i}'+ \f{\p}{\p z}R_i(z,N_{\e,i}) \vp v_{\e,i}' +m_j \vp v_{\e,i}' \Big(  \f{v_{\e,j}'-v_{\e,i}'}{\e} \Big) \exp \Big( \f{v_{\e,j}-v_{\e,i}}{\e} \Big)\\
& -\f \e 2 \vp'' p_{\e,i}-\e \vp' p_{\e,i}'-2 p_{\e,i}\vp'v_{\e,i}'.
\end{array}
$$
Let's suppose that
$$
\max_{z\in \R} (P_{\e,1}(z), P_{\e,2}(z))=P_{\e,1}(z_0),\quad \text{for $z\in B_1(\xi)$.}
$$ 
 Then, evaluating the equation on $P_{\e,1}$ at $z_0$ we obtain
$$
\e \vp(z_0) v_{\e,1}''^2 (z_0)\leq  \f{\p}{\p z}R_1(z_0,N_{\e,i}) \vp(z_0) v_{\e,1}'(z_0)  -\f \e 2 \vp'' (z_0) p_{\e,1}(z_0)-\e \vp' (z_0)p_{\e,1}'(z_0)-2 p_{\e,1} (z_0)\vp' (z_0) v_{\e,1}' (z_0).
$$
Using \fer{vp} and $0=(\vp p_{\e,1})'(z_0)=\vp'(z_0) p_{\e,1}(z_0)+\vp(z_0) p_{\e,1}'(z_0)$, we obtain 
$$
\e \vp(z_0) v_{\e,1}''^2 (z_0)\leq  \f{\p}{\p z}R_1(z,N_{\e,1})  \vp(z_0) v_{\e,1}'(z_0)+\f{3C\e}{2} |v_{\e,1}' (z_0)|^2+2C \vp(z_0)^{\f12}|v_{\e,1}' (z_0)|^3.
$$
We deduce thanks to \fer{sys-ve} and the above line that, 
$$
\begin{array}{c}
\f{\vp(z_0)}{\e} \Big(|v_{\e,1}'(z_0)|^2 + R_1(z_0,N_{\e,1}) + m_2 \exp \big( \f{v_{\e,2}(z_0)-v_{\e,1}(z_0)}{\e} \big) -m_1 \Big)^2
\leq \\
  \f{\p}{\p z}R_1(z,N_{\e,1})  \vp(z_0) v_{\e,1}'(z_0)+\f{3C\e}{2} |v_{\e,1}' (z_0)|^2+2C \vp(z_0)^{\f12}|v_{\e,1}' (z_0)|^3.
 \end{array}
$$
Since $\xi\in K$,  $R_1(z,N_{\e,1})$ and $\f{\p}{\p z}R_1(z,N_{\e,1})$ are bounded uniformly by a constant depending only on $K$. We thus deduce that there exists a constant $D=D(K)$ such that for all $\e\leq \e_0$ we have
$$
|v_{\e,1}'(z_0)|^2 \leq \f{D}{\vp(z_0)},
$$
which leads to
$$
P_{\e,1}(z_0) \leq D.
$$
Since $z_0$ was the maximum point of $P_{\e,i}$, we obtain that
$$
\vp(\xi) |v'_{\e,i}(\xi)|^2 =P_{\e,i}(\xi) \leq D.
$$
However, $\vp(\xi)=1$ and hence 
$$
|v'_{\e,i}(\xi)| \leq \sqrt{D}.
$$
It is possible to do the above computations for any $\xi\in K$ and the above bound $\sqrt{D}$, depending only on $K$, will remain unchanged. We conclude that the families $(v_{\e,i})_\e$ are   uniformly Lipschitz in $K$.\\

(v) The next step is to prove the convergence along subsequences of the families $(v_{\e,i})_\e$ to a viscosity solution of \fer{HJ-v}. Note that thanks to the previous steps we know that the families $(v_{\e,i})_\e$ are locally uniformly bounded and Lipschitz. Therefore, from the Arzela-Ascoli Theorem, they converge along subsequences to functions $v_i \in \mathrm{C}(\R)$. Moreover, we deduce from \fer{v-Har} that $v_1=v_2=v$. The fact that $v$ is a viscosity solution to \fer{HJ-v} can be derived using the method of perturbed test functions similarly to the proof of Theorem 1.1 in \cite{SM:12}.\\

(vi) We next prove \fer{W-neg}. Let's suppose in the contrary that there exists $z_0\in \R$ such that $W(z_0,N_1,N_2)>0$. Then, there exists an interval $(a_0,b_0)$ such that $z_0 \in (a_0,b_0)$ and
$W(z,N_1,N_2)>0$ for $z\in (a_0,b_0)$. We then notice that $v$ being locally uniformly Lipschitz, is differentiable almost everywhere. Let's $z_1\in (a_0,b_0)$ be a differentiability point of $v$. Then from \fer{HJ-u} we obtain that 
$$
-|v' (z_1)|^2=W(z_1,N_1,N_2),
$$
which is a contradiction with the fact that $W(z_1,N_1,N_2)>0$.\\

(vii) Finally, we prove \fer{bd-bel-N}. Note from the expression of $W(z,N_1,N_2)$ in \fer{W} and from \fer{as:r-m} that $0<\max \big(W(-\theta,0,0), W(\theta,0,0) \big)$. We assume, without loss of generality, that $0<W(-\theta,0,0)$.  Therefore, there exists an interval $(a_1,b_1)$ with $-\theta\in (a_1,b_1)$ and  $\overline\da$ such that
$$
0<W(z,N_1,N_2),\qquad \text{for all $N_1,N_2< \overline \da$, and $z\in (a_1,b_1)$}.
$$
We deduce from the above line and step (vi) that there exists $i\in \{1,2\}$ such that $N_i>\overline \da$. Without loss of generality, we suppose that $i=1$. From the fact that $(N_{\e,i})_\e$ converges to $N_i$ and from Proposition \ref{prop:Ne}-(iii) we obtain that there exists a compact set $K$ and a constant $\e_0>0$ such that
$$
 \f{\overline \da}{2} \leq    \int_K n_{\e,1}(z)dz, \qquad \text{for all $\e\leq \e_0$}.
$$
We then deduce from  \ref{prop:Ne}-(ii) that
$$
\da:=\f{\overline \da}{2C_M(K)} \leq    \int_K n_{\e,2}(z)dz \leq N_{\e,2}.
$$
This completes the proof of \fer{bd-bel-N}.
\qed

\subsection{ Convergence to the demographic equilibrium of the ESS and consequences (the proof of Theorem \ref{thm:main}) }
\label{sec:HJ-unique}

We are now ready to prove Theorem \ref{thm:main}.\\

{\bf Proof of Theorem \ref{thm:main}.} 
(i) We first prove the first part of the theorem. Note that we already proved in the previous section that as $\e\to 0$, $n_{\e,i}$ converges in the sense of measures to $n_i$ and $N_{\e,i}$ converges to $N_i$ such that $\int_\R n_i(z)dz= N_i$. Moreover, thanks to \fer{supp-n} and \fer{W-neg} we have
$$
W(z,N_1,N_2)=0,\quad \text{for $z\in \mathrm{supp} \,n_i$}\quad \text{and} ,\quad W(z,N_1,N_2)\leq 0, \quad \text{for  $z\not\in \mathrm{supp}\, n_i$.}
$$
Furthermore, one can verify using \fer{W} that  $W$ can take its maximum only at one or two points and hence the support of $n_i$ contains only one or two points. This implies indeed that $ \mathrm{supp} \,n_i$ is indeed an ESS.  We then deduce from the uniqueness of the ESS (see Theorem \ref{th:ESS}) that $n_i=n_i^*$ and $N_i=N_i^*$, for $i=1,2$ and $(n_1^*,n_2^*)$ the demographic equilibrium corresponding to the unique ESS. \\

(ii) The second part of Theorem \ref{thm:main} is immediate from it's first part and the previous subsection.\\

(iii) We first notice from part (i) that $\Omega= {\rm supp}\, n_1^*= {\rm supp}\, n_2^*$ is the unique ESS of the model. Moreover, from Corollary \ref{cor:deg}  and under condition \fer{non-deg} we obtain \fer{Wneg} and consequently
$$
 {\rm supp}\, n_1^*= {\rm supp}\, n_2^* = \{z\, |\, W(z,N_1^*,N_2^*)=0\}.
 $$
The above equalities together with \fer{supp-n} lead to \fer{aubry}. It then remains to prove that the solution of \fer{HJ}--\fer{aubry} is unique. The uniqueness of $u$ indeed derives from the fact that any negative viscosity solution of \fer{HJ} can be uniquely determined by its values at the maximum points of $W$ (\cite{PL:82}, Chapter 5). However, \fer{aubry} implies that $u=0$ at such points and hence such solution is unique.  \\

\noindent
Note indeed that restricting to a bounded domain $\mathcal O$ and following similar arguments as in \cite{PL:82}--Chapter 5, we obtain that a  viscosity solution of \fer{HJ} in the domain $\mathcal O$,  verifies
$$
u(z)=\sup\,  \{ L(y,z)+u(y)\, |\, \text{with $y$ a maximum point of $W(\cdot, N_1^*,N_2^*)$ or $y\in \p \mathcal O$}\},
$$
with
$$
\begin{array}{rl}
L(y,z)=\sup \, \{&- \int_0^T \sqrt{-W(\gamma(s),N_1^*,N_2^*) } \,ds \, |\, (T,\gamma) \text{ such that } \gamma(0)=y, \gamma(T)=z, \\
&|\f{d\gamma}{ds} |\leq 1, \, \text{a.e. in }[0,T], \; \gamma(t)\in \overline{\mathcal{O}},
 \; \forall t \in [0,T]
  \} .
  \end{array}
$$ 
Although here we have an unbounded domain, the trajectories which come from infinity do not change the value of the solution since $u$ is negative and $W$ is strictly negative for $|z|$ large enough. This allows to conclude that the solution $u$ of \fer{HJ}   is indeed determined by its values at  the maximum points of $W$. Note also that the above property is  indeed a particular case of a property from   the weak KAM theory, which is the fact that the viscosity solutions are completely determined by one value taken on each static class of the Aubry set  \cite{GC:01}.  \\
 
\section{A source and sink case}
\label{sec:sink}

In this section, we consider a particular case where there is  migration only from one habitat to the other, that is 
\beq
\label{as:sink}
m_1>0, \qquad m_2=0.
\eeq
We also assume that
\beq
\label{as:r1}
r_1-m_1>0.
\eeq
Following similar arguments to the case of migration in both directions, one can characterize the mutation, selection and migration equilibria. However, since the migration is only in one direction, we should study the equilibria in the two habitats separately. \\

\noindent
Note that since $m_2=0$, there is no influence of the second habitat on the first habitat. One can indeed compute explicitly $n_{\e,1}$:
\beq
\label{n1-sink}
n_{\e,1}(z)=\f{g_1^{\f 14}N_{\e,1}}{\sqrt{2\pi\e}}\, \exp \Big(-\f{\sqrt{g_1}(z+\theta)^2}{2\e} \Big),\qquad N_{\e,1}=\f{r_1-m_1-\e\sqrt{g_1}}{\kappa_1}.
\eeq
Note that as $\e\to 0$, $n_{\e,1}$ converges in the sense of measures to $n_1^{M*}$ with
$$
n_1^{M*}(z)=N_1^{M*} \da(z+\theta),\qquad N_1^{M*} = \f  {r_1-m_1}{\kappa_1}.
$$
Here, $\{-\theta\}$ is indeed the unique ESS in the first habitat and $n_1^*$ corresponds to the demographic equilibrium at the ESS. \\

\noindent
In the second habitat however, there is an influence of the population coming from the first habitat. The natural quantity that appears in this case as the effective fitness in the second habitat is still the principal eigenvalue of \fer{efitness} which is, in this case, given by
$$
\begin{array}{rl}
W(z,N_2)&=\max( r_1-g_1(z+\theta)^2-\kappa_1N_1^{M*}-m_1, r_2-g_2(z-\theta)^2-\kappa_2N_2)\\
& =\max(  -g_1(z+\theta)^2 , r_2-g_2(z-\theta)^2-\kappa_2N_2).
\end{array}
$$
Then one can introduce the notion of the ESS for this habitat similarly to Section \ref{sec:ad}.

\subsection{The results in the adaptive dynamics framework }

We can indeed always identify the unique ESS:

\begin{theorem}
\label{th:ESS2}
Assume \fer{as:sink}--\fer{as:r1}. In each patch there exists a unique ESS. In patch $1$ the ESS is always monomorphic and it is given by $\{-\theta\}$ with the following demographic equilibrium:
\beq
\label{asdim-sink}
n_1^{M*}= N_1^{M*} \, \da(z+\theta),\quad  N_1^{M*}=\f{r_{1}-m_1}{\kappa_1}.
\eeq

In patch $2$ there are two possibilities:\\
(i)  the ESS is 
dimorphic if and only if
\beq
\label{con-dim-source}
\f{m_1(r_{1}-m_1)}{\kappa_1} < \f{4g_2\theta^2 r_{2}}{\kappa_2}.
\eeq
The dimorphic ESS is given by $\{-\theta,\theta\}$ with the following demographic equilibrium:
$$
n_2^{D*}=\alpha \da(z+\theta)+\beta \da(z-\theta),\quad N_2^{D*}=\alpha+\beta= \f{r_{2}}{\kappa_2},
\quad
\alpha =\f{m_1(r_{1}-m_1)}{4g_2\theta^2 \kappa_1},\quad \beta =\f{r_{2}}{\kappa_2}-\f{m_1(r_{1}-m_1)}{4g_2\theta^2 \kappa_1}.
$$

(ii) If condition \fer{con-dim-source} is not satisfied then the ESS in the second patch is monomorphic. The ESS is given by $\{-\theta\}$  with the following demographic equilibrium:  
$$
n_2^{M*}=  N_2^{*}\, \da(z+\theta),
\quad
N_2^{M*} = \f{1}{2\kappa_2}\Big( r_{2}-4g_2\theta^2 
+\sqrt{ (r_{2}-4g_2\theta^2)^2  +4 \f{\kappa_2}{\kappa_1}m_1  (r_{1}-m_1) }
\Big).
$$

\end{theorem}

\noindent
The proof of the above theorem is not difficult and is left to the interested reader. 

\subsection{The computation of the zero order term $u_2$}

We then proceed with the method presented in the introduction to characterize the evolutionary equilibrium $n_{\e,2}(z)$. To this end, we first identify the zero order term $u_2$ (introduced in \fer{WKB}--\fer{ap-ue}):

\begin{theorem}\label{thm:sink}
 Assume \fer{as:sink}--\fer{as:r1}.
\\
(i) As $\e \to 0$, $(n_{\e,1}, n_{\e,2})$ converges to $(n_1^{M*},n_2^{*})$, the demographic equilibrium of the unique   ESS of the metapopulation, given by Theorem \ref{th:ESS2}. Moreover, as $\e \to 0$, $(N_{\e,1}, N_{\e,2})$ converges to $(N_1^{M*},N_2^*)$, the total populations in patch $1$ and $2$ corresponding to this demographic equilibrium.
\\
(ii)  As $\e\to 0$,  $(u_{\e,2})_\e$  converges  locally uniformly in $\R$ to   $u_1(z) =- \f{\sqrt{g_1}}{2}(z+\theta)^2$. 
As $\e\to 0$,  $(u_{\e,2})_\e$ converges along subsequences and locally uniformly in $\R$ to a  function $u_2\in \mathrm{C}(\R)$ which satisfies 
\begin{equation}
\label{HJ-source-2}
 -|  u_2'|^2\leq  \max(R_1(z,N_1^{M*})-m_1, R_2(z,N_2^*)),   \quad
  -| u_2'|^2 \geq  R_2(z,N_2^*),
   \quad u_1(z) \leq u_2(z),
   \quad \max_{z\in \R}u_2(z)=0,
\end{equation}
where the first two inequalities are in the viscosity sense. Moreover, we have the following condition on  the zero level set of $u_2$:
\beq
\label{n2supp}
  {\rm supp}\, n_2^*\subset \{z\, |\, u_2(z)=0 \} \subset \{z\, |\, \max(R_1(z,N_1^{M*})-m_1, R_2(z,N_2^*))=0\}.
\eeq
 \end{theorem}
 
 \noindent
{\bf Proof.} The proof of Theorem \ref{thm:sink} is close to the proof of Theorem \ref{thm:main}-(i) and (ii). We only provide the steps of the proof and discuss the main differences.\\
(i) We first notice that the convergence of $(n_{\e,1})_\e$, $(N_{\e,1})_\e$ and $(u_{\e,1})_\e$ is trivial from \fer{n1-sink}.\\

\noindent
(ii) Following similar arguments as in Proposition \ref{prop:Ne}--(i) and (iii) we find that $N_{\e,2}$ is bounded from above and that $n_{\e,2}$ has small mass at infinity. Hence, as $\e\to 0$ and along subsequences, respectively $(n_{\e,2})_\e$ and $(N_{\e,2})_\e$ converges to $n_2$ and $N_2$ with $N_2=\int n_2(z)dz$.\\

\noindent
(iii) Note that since $m_2=0$, \fer{Harnack} does not hold anymore but a weaker version of it still holds true. We can indeed obtain, following similar arguments and still referring to \cite{JB.MS:04}, Theorem 8.2, that for any compact set $K\subset \R$, there exists indeed a constant $C_M=C_M(K)$ such that, for all $\e\leq 1$,we  have 
\beq
\label{Harnack-sink}
 n_{\e,1}(x) \leq C_M n_{\e,2}(y), \qquad \text{for}\quad   |x-y|\leq \e.
\eeq

\noindent
(iv) We deduce from \fer{Harnack-sink} and the fact that $n_{\e,1}$ has small mass at infinity, that there exists $\e_0$ such that, for all $\e\leq \e_0$, $N_{\e,2}$ is uniformly bounded from below by a positive constant.\\

\noindent
(v) Following similar arguments as in the proof of Proposition \ref{prop:ve} we obtain that there exists $\e_0>0$, such that for all $\e\leq \e_0$, $(u_{\e,2})_\e$ is locally uniformly bounded and Lipschitz. Therefore, as $\e\to 0$ and along subsequences, $(u_{\e,2})_\e$ converges to a  function $u_2\in \mathrm{C}(\R)$ such that $\max_{z\in \R}u_2(z)=0$. Moreover, from \fer{Harnack-sink} we obtain that $u_1(z) \leq u_2(z)$, for all $z\in\R$.\\

\noindent
(vi) Note that $u_{\e,2}$ solves the following equation
\beq
\label{ue2}
-\e u_{\e,2}''=|u_{\e,2}'|^2+R_2(z,N_{\e,2}) + m_1 \exp \Big( \f{u_{\e,1}-u_{\e,2}}{\e} \Big).
\eeq
Passing to the limit as $\e\to 0$ and using the fact that the last term above is positive we obtain that
$$
-|u_{2}'|^2 \geq R_2(z,N_2),
$$
in the viscosity sense.

(vii) Next, we prove that 
$$
 -|  u_2'|^2\leq\max(R_1(z,N_1^{M*})-m_1, R_2(z,N_2)).
$$
To this end, we consider two cases. Let's first suppose that $z_0$ is such that $u_2(z_0)=u_1(z_0)$. Moreover, let $\vp$ be a smooth test function such that $u_2-\vp$ has a local maximum    at $z_0$. Then, since $u_1(z) \leq u_2(z)$, $u_1-\vp$ has also a local maximum  at $z_0$ and hence 
$$
-|\vp'(z_0) |^2\leq R_1(z,N_1^{M*})-m_1\leq \max(R_1(z,N_1^{M*} )-m_1, R_2(z,N_2)).
$$
Next we assume that $u_1(z_0) < u_2(z_0)$. In this case, as $\e\to 0$, the last term in \fer{ue2} tends to $0$ at $z_0$ and hence
$$
-|u_2'(z_0) |^2\leq R_2(z,N_2) \leq \max(R_1(z,N_1^{M*} )-m_1, R_2(z,N_2)),
$$
in the viscosity sense.\\

(viii) We then prove \fer{n2supp}. The fact that $ {\rm supp}\, n_2\subset \{z\, |\, u_2(z)=0 \}$ is immediate from \fer{WKB}. To prove the second property, we first notice that, considering $0$ as a test function, $ -|  u_2'|^2\leq\max(R_1(z,N_1)-m_1, R_2(z,N_2))$ implies that 
$$
0 \leq \max(R_1(z,N_1^{M*} )-m_1, R_2(z,N_2)), \qquad \text{in $ \{z\, |\, u_2(z)=0 \}$.}
$$
Moreover, $-|u_{2}'|^2 \geq R_2(z,N_2)$, implies that $R_2(z,N_2)\leq 0$. We also know that $R_1(z,N_1^{M*} )-m_1\leq 0$.  Hence,  \fer{n2supp}.

(ix) Finally, we deduce from the previous step that 
$$
W(z,N_2)\leq 0, \text{ in $\R$}, \qquad W(z,N_2)=0, \text{ for $z\in  {\rm supp}\, n_2$}.
$$
This means that ${\rm supp}\, n_2$ is an ESS and hence, thanks to Theroem \ref{th:ESS2}, we obtain   that $n_2=n_2^*$ and $N_2=N_2^*$, where $n_2^*$ and $N_2^*$ are given by Theorem \ref{th:ESS2}. We then deduce in particular that the whole sequences $(n_{\e,2})_\e$ and $(N_{\e,2})_\e$ converge respectively to $n_2^*$ and $N_2^*$.

\qed

Theorem \ref{thm:sink} allows us to identify $u$ in a neighborhood of the ESS points:

\begin{prop}
\label{prop:expu-sink}
(i) There exists a connected and open set $\mathcal O_{\mathrm I}\subset \R$, with $-\theta\in \mathcal O_{\mathrm I}$, such that
$$
u_2(z)=-\f{\sqrt{g_1}}{2}(z+\theta)^2.
$$
(ii) Assume that \fer{con-dim-source} holds. Then, there exists a connected and open set $\mathcal O_{\mathrm{II}}\subset \R$, with $\theta\in \mathcal O_{\mathrm{II}}$, such that
$$
u_2(z)=-\f{\sqrt{g_2}}{2}(z-\theta)^2.
$$
(iii) Assume that 
\beq
\label{mrm}
\f{4g_2\theta^2 r_{2}}{\kappa_2}<\f{m_1(r_{1}-m_1)}{\kappa_1} .
\eeq
Then $u_2(\theta)<0$.

\end{prop}
Note that when $\f{m_1(r_{1}-m_1)}{\kappa_1}  = \f{4g_2\theta^2 r_{2}}{\kappa_2}$ we don't know the value of $u_2(\theta)$. In particular, it can vanish. This is why we cannot provide an approximation of $n_{\e,2}$ in this degenerate case. 
\\

\noindent
{\bf Proof of Proposition \ref{prop:expu-sink}.} (i) Note that using similar arguments as in the proof of Theorem \ref{thm:main}-(iii), where we used properties from the weak KAM theory, and using 
$$
-|u_2'|^2 (z) \leq W(z;N_2^*)=\max(R_1(z,N_1^{M*} )-m_1, R_2(z,N_2^*))\leq 0, \quad u_2(z)\leq 0, 
$$
which holds thanks to  \fer{HJ-source-2},
we obtain that
\beq
\label{u2leq}
 u_2(z)\leq \max \Big( - |\int_{{\theta}}^{ z} \sqrt{- W(x;   N_2^{*})} dx|
  , - |\int_{-\theta}^{z} \sqrt{- W(x;  N_2^{*})} dx |\Big).
\eeq
From the above inequality it is immediate that there exists a connected and open set $\mathcal O_{\mathrm I}\subset \R$, with $-\theta\in \mathcal O_{\mathrm I}$, such that
$$
u_2(z)\leq u_1(z)=-\f{\sqrt{g_1}}{2}(z+\theta)^2.
$$
Combining this with the third property in \fer{HJ-source-2} we deduce the first claim of Proposition \ref{prop:expu-sink}.\\

\noindent
(ii) Note that under condition  \fer{con-dim-source} the ESS is dimorphic and that $\mathrm{supp} \, n_2^{D*}=\{-\theta,\theta\}$. Therefore, we deduce thanks to \fer{n2supp} that
$u_2(\theta)=0$. This property combined with the second property in \fer{HJ-source-2} implies that
$$
u_2(z) \geq -\f{\sqrt{g_2}}{2}(z-\theta)^2.
$$
The second claim of the theorem then follows from \fer{u2leq}.\\

\noindent
(iii) Finally, we prove the third claim of the theorem. To this end, we assume    that \fer{mrm} holds, and hence the ESS in the second patch is monomorphic and given by $\{-\theta\}$, but $u_2(\theta)=0$. Note that similarly, to the case of migration in both directions, $u_2$ is a semiconvex function. Therefore it is differentiable at its maximum points and in particular at $\theta$. Hence, the first claim of \fer{HJ-source-2} implies that
$$
0 \leq \max(R_1(\theta,N_1^{M*} )-m_1, R_2(\theta,N_2^{M*})) .
$$
However, this is in contradiction with  \fer{mrm}.   
\qed

\subsection{Next order terms}

In this subsection we compute the next order terms in the approximation of $u_{\e,i}$ and $N_{\e,i}$:
$$
u_{\e,i}=u_i+\e v_i+\e^2w_i+o(\e^2),\qquad N_{\e,i}=N_i^*+\e K_i+O(\e^2).
$$
We first notice that, thanks to \fer{n1-sink} we already know explicitly $u_{\e,1}$ and $N_{\e,1}$:
$$
u_{\e,1}= -\f{\sqrt{g_1}(z+\theta^2)}{2}+\e \log \Big(   g_1^{\f14} \big(N_1^{M*}-\e\f{\sqrt{g_1}}{\kappa_1} \big)\Big),\quad N_{\e,1}=\f{r_1-m_1-\e\sqrt{g_1}}{\kappa_1},
$$ 
and hence
\beq
\label{v1-sink}
v_1\equiv \log \big( g_1^{\f 14}N_1^{M*} \big), \qquad w_1\equiv  -\f{\sqrt{g_1}}{\kappa_1 N_1^{M*}},\qquad K_1=- \f{\sqrt{g_1}}{\kappa_1} .
\eeq

\noindent
We next compute $v_2$ and $w_2$ around the ESS points. 
\noindent
We only present the method to compute $v_2$ and $w_2$ around $-\theta$, in  the case where 
$$
\f{m_1(r_{1}-m_1)}{\kappa_1} > \f{4g_2\theta^2 r_{2}}{\kappa_2},
$$
so that the ESS is monomorphic and is given by $\{-\theta\}$. The dimorphic case, where  \fer{con-dim-source} is satisfies,   can be analyzed  following similar arguments. We recall that in the degenerate case where $\f{m_1(r_{1}-m_1)}{\kappa_1} = \f{4g_2\theta^2 r_{2}}{\kappa_2},$ we don't provide an approximation of $n_{\e,2}$. \\

\noindent
To compute $v_2$, we keep the zero order terms in \fer{ue2} in $\mathcal O_I$ and using \fer{v1-sink} we obtain
$$
v_2(z)=\log \Big(\f{m_1g_1^{\f 14}N_1^{M*}}{-g_1(z+\theta)^2+g_2(z-\theta)^2-r_2+\kappa_2 N_2^{M*}} \Big),\quad \text{for $z\in \mathcal O_I$}.
$$
Similarly to Section \ref{sec:vw} we write a Taylor expansion for $v_2$ around $-\theta$:
$$
v_2(z)=v_2(-\theta)+D_2(z+\theta)+ E_2 (z+\theta)^2+O(z+\theta)^3,\quad \text{with } v_2(-\theta)=\log (g_1^{\f14} N_2^{M*}),
$$
and we define $w_2(-\theta)=F_2$. Note that $D_2$ and $E_2$ are known thanks to the explicit computation of $v_2(z)$ given above. 
Similarly to Section \ref{sec:vw}, keeping the first order terms in $\f{1}{\sqrt{2\pi\e}}\int_{\mathcal I} \exp \big( \f{u_{\e,2}(z)}{\e} \big)dz$ we obtain that 
\beq
\label{K2-sink}
K_2=N_2^{M*}\big( \f{E_2+0.5D_2^2}{\sqrt{g_1}}+F_2 \big).
\eeq
Moreover, keeping the first order terms in  \fer{ue2} in $\mathcal O_I$ we obtain that
\beq
\label{w2-sink}
\sqrt g_1= -2\sqrt{g_1} (z+\theta) v_2'-\kappa_2K_2 + m_1 \f{N_1^{M*}}{N_2^{M*}}(w_1-w_2).
\eeq
We evaluate the above equation at $-\theta$ to obtain
$$
F_2=-\f{\sqrt {g_1}}{\kappa_1 N_1^{M*}}-\f{N_2^{M*}}{m_1N_1^{M*}} \big( \sqrt{g_1}+\kappa_2K_2 \big).
$$
One can then compute $K_2$ and $F_2$ combining the above equation with \fer{K2-sink}. Note finally that, once $K_2$ is known, one can compute $w_2$ in $\mathcal I$ thanks to \fer{w2-sink}. 

%
%

\section*{Acknowledgements}
 The author is immensely thankful to Sylvain Gandon for introducing the biological motivations to her and for determining the biological directions of this work. This article is a part of a project that we derived together with the objective of  introducing the presented method to the biological community.
The author  is also grateful for partial funding from the European Research Council (ERC) under the European Union's Horizon 2020 research and innovation programme (grant agreement No 639638), held by Vincent Calvez, and from the french ANR projects KIBORD ANR-13-BS01-0004 and MODEVOL ANR-13-JS01-0009.


\begin{thebibliography}{}

\end{thebibliography}


\begin{thebibliography}{10}

\bibitem{JB.MS:04}
J.~Busca and B.~Sirakov.
\newblock Harnack type estimates for nonlinear elliptic systems and
  applications.
\newblock {\em Ann. Inst. H. Poincar\'e Anal. Non Lin\'eaire}, 21:543--590,
  2004.

\bibitem{GC:01}
G.~Contreras.
\newblock Action potential and weak kam solutions.
\newblock {\em Calc. Var. Partial Differential Equations}, 13(4):427--458,
  2001.

\bibitem{TD:00}
T.~Day.
\newblock Competition and the effect of spatial resource heterogeneity on
  evolutionary diversification.
\newblock {\em The American Naturalist}, 155(6):790--803, 2000.

\bibitem{FD.OR.SG:13}
F.~D\'ebarre, O.~Ronce, and S.~Gandon.
\newblock Quantifying the effects of migration and mutation on adaptation and
  demography in spatially heterogeneous environments.
\newblock {\em Journal of Evolutionary Biology}, 26:1185--1202, 2013.

\bibitem{OD.PJ.SM.BP:05}
O.~Diekmann, P.-E. Jabin, S.~Mischler, and B.~Perthame.
\newblock The dynamics of adaptation: an illuminating example and a
  {H}amilton-{J}acobi approach.
\newblock {\em Th. Pop. Biol.}, 67(4):257--271, 2005.

\bibitem{LE.PS:89}
L.~C. Evans and P.~E. Souganidis.
\newblock A {PDE} approach to geometric optics for certain semilinear parabolic
  equations.
\newblock {\em Indiana Univ. Math. J.}, 38(1):141--172, 1989.

\bibitem{CF.SM.EM:12}
C.~Fabre, S.~M\'el\'eard, E.~Porcher, C.~Teplitsky, and Robert A.
\newblock Evolution of a structured population in a heterogeneous environment.
\newblock {\em Preprint}.

\bibitem{AF:16}
A.~Fathi.
\newblock {\em Weak Kam Theorem in Lagrangian Dynamics}.
\newblock Number~88 in Cambridge Studies in Advanced Mathematics. Cambridge
  Univ Pr, 2016.

\bibitem{MF:86}
M.~Freidlin.
\newblock Geometric optics approach to reaction-diffusion equations.
\newblock {\em SIAM J. Appl. Math.}, 46:222--232, 1986.

\bibitem{SG.SM:16}
S.~Gandon and S.~Mirrahimi.
\newblock A {H}amilton-{J}acobi method to describe the evolutionary equilibria
  in heterogeneous environments and with non-vanishing effects of mutations.
\newblock {\em Preprint}.

\bibitem{AH.TD.ET:01}
A.~Hendry, T.~Day, and E.~B. Taylor.
\newblock Population mixing and the adaptive divergence of quantitative traits
  in discrete populations: a theoretical framework for empirical tests.
\newblock {\em Evolution}, 55(3):459--466, 2001.

\bibitem{PL:82}
P.~L. Lions.
\newblock {\em Generalized solutions of {H}amilton-{J}acobi equations},
  volume~69 of {\em Research Notes in Mathematics}.
\newblock Pitman Advanced Publishing Program, Boston, 1982.

\bibitem{GM.IC.SG:97}
G.~Mesz{\'e}na, I.~Czibula, and S.~Geritz.
\newblock Adaptive dynamics in a 2-patch environment: a toy model for
  allopatric and parapatric speciation.
\newblock {\em Journal of Biological Systems}, 5(02):265--284, 1997.

\bibitem{SM:11}
S.~Mirrahimi.
\newblock {\em Concentration phenomena in PDEs from biology}.
\newblock PhD thesis, Univeristy of Pierre et Marie Curie (Paris 6), 2011.

\bibitem{SM:12}
S.~Mirrahimi.
\newblock Migration and adaptation of a population between patches.
\newblock {\em Discrete and Continuous Dynamical Systems - Series B (DCDS-B)},
  18(3):753--768, 2013.

\bibitem{SG.SM:17}
S.~Mirrahimi and S.~Gandon.
\newblock The equilibrium between selection, mutation and migration in
  spatially heterogeneous environments.
\newblock {\em In preparation}.

\bibitem{SM.JR:16}
S.~Mirrahimi and Roquejoffre J.-M.
\newblock A class of {H}amilton-{J}acobi equations with constraint: Uniqueness
  and constructive approach.
\newblock {\em Journal of Differential Equations}, 260(5):4717--4738, 2016.

\bibitem{SM.JR:15-1}
S.~Mirrahimi and J.-M. Roquejoffre.
\newblock Uniqueness in a class of {H}amilton-{J}acobi equations with
  constraints.
\newblock {\em C. R. Math. Acad. Sci. Paris}, 353:489--494, 2015.

\bibitem{GB.BP:08}
B.~Perthame and G.~Barles.
\newblock Dirac concentrations in {L}otka-{V}olterra parabolic {PDE}s.
\newblock {\em Indiana Univ. Math. J.}, 57(7):3275--3301, 2008.

\bibitem{Rice-book}
S.~H. Rice.
\newblock {\em Evolutionary theory: mathematical and conceptual foundations}.
\newblock Sinauer Associates, Inc., 2004.

\bibitem{OR.MK:01}
O.~Ronce and M.~Kirkpatrick.
\newblock When sources become sinks: migration meltdown in heterogeneous
  habitats.
\newblock {\em Evolution}, 55(8):1520--1531, 2001.

\bibitem{SY.FG:09}
S.~Yeaman and F.~Guillaume.
\newblock Predicting adaptation under migration load: the role of genetic skew.
\newblock {\em Evolution}, 63(11):2926--2938, 2009.

\end{thebibliography}
\bibliographystyle{plain}

\end{document}